\documentclass[11pt,leqno]{article}
\usepackage{amsthm,amsfonts,amssymb,epsfig,graphics,amsmath}
\usepackage[latin1]{inputenc}\relax

\newcommand{\sgn}{{\text{\rm sgn}}}

\newcommand{\calB}{{\cal{B}}}
\newcommand{\calE}{{\cal{E}}}

\newcommand{\R}{\Re}

\newcommand{\CC}{{\mathbb C}}
\newcommand{\mA}{{\mathbb A}}

\newcommand\cE{{\cal  E}}

\newcommand\adots{\mathinner{\mkern2mu\raise1pt\hbox{.}
\mkern3mu\raise4pt\hbox{.}\mkern1mu\raise7pt\hbox{.}}}

\newtheorem{theo}{Theorem}[section]
\newtheorem{prop}[theo]{Proposition}
\newtheorem{cor}[theo]{Corollary}
\newtheorem{lem}[theo]{Lemma}

\newtheorem{ass}[theo]{Assumption}

\newtheorem{rem}[theo]{Remark}
\newtheorem{rems}[theo]{Remarks}

\numberwithin{equation}{section}
\newcommand\be{\begin{equation}}
\newcommand\ee{\end{equation}}
\newcommand\bp{\begin{pmatrix}}
\newcommand\ep{\end{pmatrix}}
\newcommand\ba{\begin{aligned}}
\newcommand\ea{\end{aligned} }

\newcommand\bea{\begin{array}}
\newcommand\ena{\end{array}}

\begin{document}

\title{\bf Existence and stability of steady states
of a reaction convection diffusion equation modeling
microtubule formation}
\author{\sc \small
Shantia Yarahmadian\thanks{Mississippi State University, MS 39762;
syarahmadian@math.msstate.edu:
Research of S.Y. was partially supported
under NSF grants number DMS-0070765 and DMS-0300487} and,
Blake Barker\thanks{ Brigham Young University, Provo, UT 84602;
bhbarker@gmail.com: Research of B.B. was
partially supported under NSF grants number DMS-0607721 and DMS-0300487},
Kevin Zumbrun\thanks{Indiana University, Bloomington, IN 47405;
kzumbrun@indiana.edu:
Research of K.Z. was partially supported
under NSF grants number DMS-0070765 and DMS-0300487.}
Sidney L. Shaw\thanks{Indiana University, Bloomington, IN 47405, sishaw@indiana.edu:
Research of S.S. was partially supported by the Indiana Metacyte Institute at Indiana University, founded in part
through a major grant from the Lily Endowment, INC.
}
}
\maketitle
\begin{abstract}
We generalize the Dogterom-Leibler model for microtubule dynamics \cite{DL} to the case
where  the rates of elongation as well as the lifetimes of the elongating
and shortening phases are a function of GTP-tubulin concentration.
We study also the effect of nucleation rate in the form of a damping
term which leads to new steady-states.
For this model, we study existence and stability of steady states
satisfying the boundary conditions at $x=0$.
Our stability analysis introduces numerical and analytical Evans function
computations as a
new mathematical tool in the study of microtubule dynamics.
\end{abstract}

\tableofcontents
\paragraph{General Theory}
We analyze some mathematical aspects of the
phenomenon of dynamic instability of microtubules. Our theoretical model includes a
simplified, semi-infinite geometry, in which infinitely
rigid microtubules grow perpendicularly to a nucleating planar
surface (See Figure 2).
Dogterom and Leibler  studied this theoretical model with this semi-infinite geometry, neglecting concentration variations in the process of assembly and disassembly of microtubules \cite{DL}. In this model infinitely rigid microtubules are growing perpendicular to a nucleation planar surface (Fig. 2) and  randomly switching between the growth state ($+$) and shrinkage state($-$) or polymerization and depolymerization.\\

In this model each microtubule switches randomly between
the assembly state $+$, in which it grows with the
average speed $v^+$ proportional to the local monomer density
$c$, and the disassembly state $-$, in which it shrinks
with the average speed $v^-$. The frequencies of the transitions
between the two states, $f_+^-$ (of the "catastrophies"
from $+$ state to $-$ state) and $f_-^+$ (the "rescues" from
$-$ state to $+$ state), determine, together with $v^+$, $v^-$ determine
the behavior of the
model.
We summarize now the main results obtained
through Evan function and time-evolution simulations.
The main result is the existence and stability of steady state
solutions for a convection  diffusion equation modeling the
above-described process.
The novelty in this study is the analysis of the general case when the dynamics parameter are depending on local concentration.
The mathematical tools introduced to handle the more complicated
resulting equations should be of general use.
\section{Introduction}\label{intro}

Microtubules are natural polymers found inside of living eukaryotic cells.  Each polymer subunit is an obligate dimer of proteins synthesized from the alpha and beta tubulin genes.  When concentrated above a threshold level, tubulin dimers bind to each other forming a beta lattice of subunits that bind head to tail and side-to-side.  Under appropriate conditions, the inherent curvature of the lattice sheet produces a tube, typically 13 subunits in cross-section, which can extend at the tube ends by polymerization to many thousands of subunits \cite{MK}.  The resultant MT polymer serves as a semi-rigid structural element in the cell that is critical for intracellular distribution networks, chromosome segregation before cell division, and neuronal activity.
%
MTs within cells are typically initiated from a complex of proteins that form a nucleation template.  Polymerization proceeds at a rate that is dependent upon the concentration of free tubulin subunits.  Over time, a remarkable phenomenon is observed both in cellular MTs and with purified tubulin.  The MTs undergo stochastic switching between states of polymerization and depolymerization, a process termed `dynamic instability' \cite{SGCOH}.  Switching frequencies are modulated during certain cellular transitions to alter the length distribution and density of the polymer array.  The mechanisms governing dynamic instability are still an active subject of both experimental and theoretical investigation \cite{KM}.

In general terms, the dynamic instability arises due to an energy dependent mechano-chemical transition that occurs in the beta-tubulin half of the tubulin dimer after polymerization.
MT polymerization proceeds in a head to tail fashion from the nucleation template leaving the beta-tubulin face exposed at the microtubule end.  Free tubulin dimers rapidly bind guanosine triphosphate (GTP), a small molecule used as a convertible form of stored chemical energy in the cell.  If the exposed beta-tubulin on the MT end is GTP-bound, it has a relatively high affinity for binding the alpha-tubulin face of any free subunits in the cytoplasm.  Once bound, the beta-tubulin undergoes a structural change increasing the probability of hydrolyzing the associated GTP to guanosine diphosphate (GDP).  The binding affinity of GDP-bound subunits for other tubulin dimers is relatively low.  Therefore, GTP hydrolysis serves as a switch for changing the binding affinity of the subunit. If the rate of GTP hydrolysis in the MT outpaces binding of new GTP-bound dimers, the MT will lose its theorized cap of GTP-bound dimers and transit from a state of high affinity dimer binding to very low affinity.  Over an important range of free tubulin concentrations, this loss of the GTP-cap will result in a switch from polymerization to depolymerization, termed `catastrophe.'  Should a GTP-bound dimer manage to bind the depolymerizing end, the polymerization can be recovered, a process termed `rescue.'  The observed stochastic switching between states of growth and shortening is attributed to the statistical variance in the processes associated with dimer binding and GTP hydrolysis.

A closed form kinetic description of microtubule dynamics has not been successfully rendered to date.  Several analytical models have been developed to describe how the major factors leading to dynamic instability (i.e. growth and shortening velocities together with catastrophe and rescue frequencies) will produce a steady state system of polymers under various conditions.
The stochastic nature of the switching between polymerization states complicates the construction of determinative equations that can be numerically solved for properties such as MT length distribution.  In this work, we advance the analysis of MT dynamics by extending the analytical methods applied to the most prominent model for dynamics instability in an in vitro context \cite{DL}.

\section{Preliminaries}
\subsection{Dichotomous Markov Noise}
The dichotomous Markov noise (DMN) $v(t)$, is a simple two-valued stochastic process such that the state space of the random variable
consists only two values $\{\pm v^{\pm}\}$
with constant transition rates $f_-^+$ and $f_+^-$. This means that the waiting-times in the two states are exponentially-distributed stochastic variables. The switches of $v(t)$ are Poisson process with probabilities $\mathfrak{p}^-$ and $\mathfrak{p}^+$ \cite{IB}. We can describe these probabilities by a first-order kinetic equation:

\begin{figure}[h!]
\centering
\includegraphics[scale=.5]{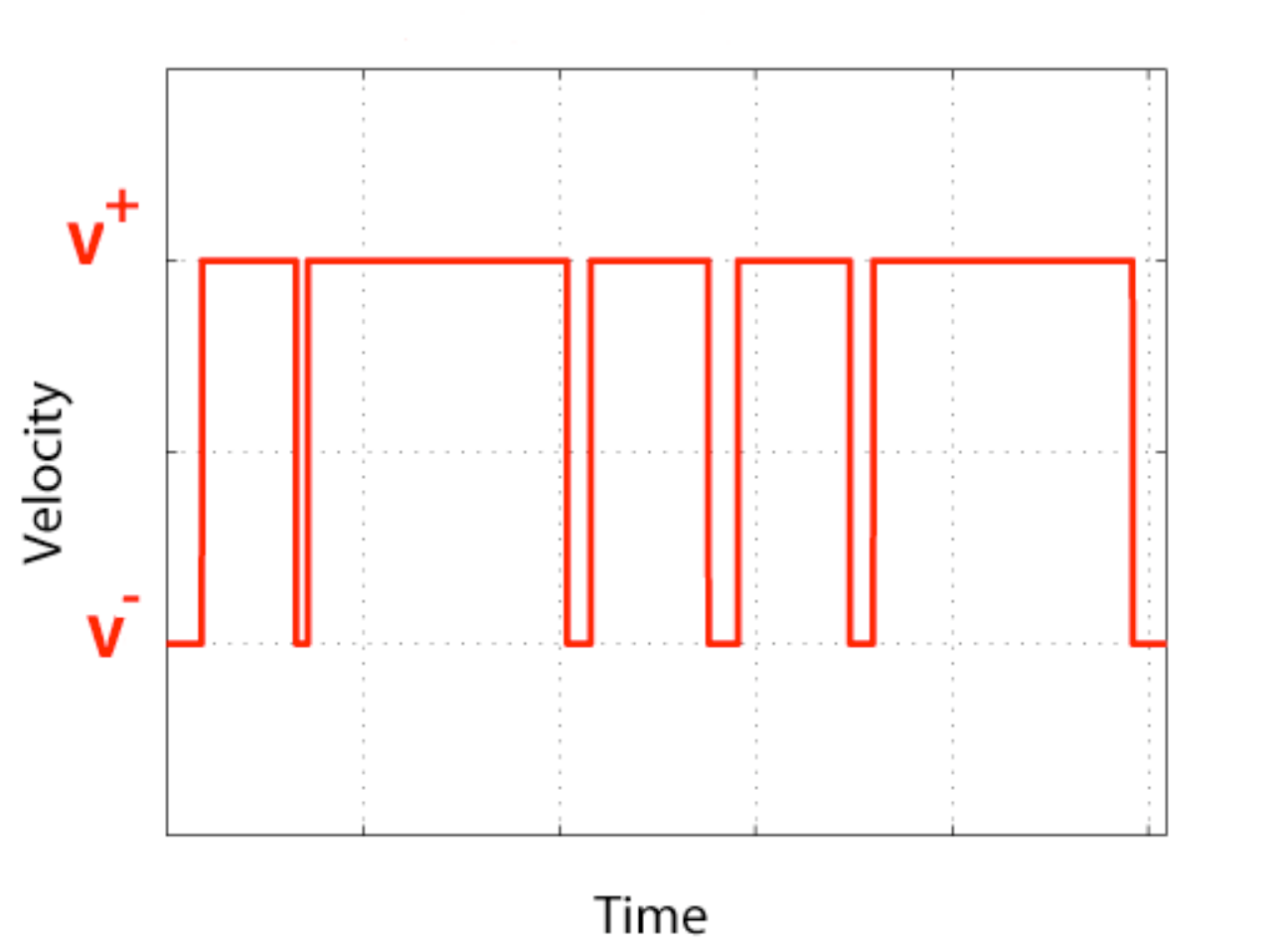}
\caption{Dichotomous Markov Noise.}
\end{figure}

\be
\frac{d}{dt}\bp \mathfrak{p}^+ \\ \mathfrak{p}^- \ep=\bp -f_+^- && f_-^+ \\  f_+^-  && -f_-^+ \ep
\bp \mathfrak{p}^+ \\ \mathfrak{p}^- \ep
\ee
\be
\mathfrak{p}^+(t) =\frac{f_-^+}{\mathcal{F}}+ \Big(\mathfrak{p}^+(0)-\frac{f_-^+}{\mathcal{F}}\Big)e^{-\mathcal{F}t}
\ee
\be
\mathfrak{p}^-(t) =\frac{f_+^-}{\mathcal{F}}+ \Big(\mathfrak{p}^-(0)-\frac{f_+^-}{\mathcal{F}}\Big)e^{-\mathcal{F}t}
\ee
where $\mathfrak{p}^+(t)+\mathfrak{p}^-(t)=1$ for all times, and $\tau_c=\frac{1}{\mathcal{F}}$=$\frac{1}{f_+^-+f_-^+}$ shows the mean time between switches of $v(t)$ or the velocity relaxation time.
The instantaneous average velocity is defined as:
\be
\overline{v}(t)=v^+\mathfrak{p^+}(t)-v^-\mathfrak{p^-(t)}=\overline{v}(\infty)+\Big(\overline{v}(0)-\overline{v}(\infty)\Big)e^{-{\mathcal{F}}t}
\ee
Where :
\be
\overline{v}(0)=v^+\mathfrak{p}^+(0)-v^-\mathfrak{p}^-(0)
\ee
and
\be
\overline{v}(\infty)=\frac{v^+f_-^+-v^-f_+^-}{\mathcal{F}}
\ee
This noise is a time-homogeneous Markov process and is therefore completely characterized by the following transition probabilty:
$$P_{ij}(t)=Pr(v(t)=i|v(0)=j), \qquad i,j \in \{\pm v^{\pm}\} $$
The temporal evolution of the noise is given by the Kolmogorov forward equation, in the physical literature known as master equation \cite{IB}:
\begin{equation}
\frac{d}{dt}\bp P_{-j}(t)\\ P_{+j}(t)\ep=\bp -f_+^- && f_-^+ \\  f_+^-  && -f_-^+ \ep\bp P_{-j}(t)\\ P_{+j}(t)\ep
\end{equation}
where $k_+$ and $k_-$ are the mean frequencies of passage from $A_+$ to $-A_-$
This system can be described by its transition matrix:
\begin{equation}
\bp
P_{--}(t)&P_{-+}(t)\\
P_{+-}(t)&P_{++}(t)
\ep=\tau_c \bp f_+^-+f_-^+e^{-\frac{t}{\tau_c}}&&f_+^-(1-e^{-\frac{t}{\tau_c}})\\
f_-^+(1-e^{-\frac{t}{\tau_c}})&&f_-^++f_+^-e^{-\frac{t}{\tau_c}} \ep
\end{equation}
 In the stationary case we have:
\be \label{stationary}
Pr(v=v^+)=f_-^+\tau_c \quad Pr(v=-v^-)=f_+^-\tau_c
\ee

We also assume that the external noise is a stationary random process and hence the DMN has to be started with
\eqref{stationary} as initial condition.
The corresponding stationary mean-value is :
\be
<v(t)>=V=(v^+f_-^+-v^-f_+^-)\tau_c
\ee
\medskip
\subsection{Master Equations}
Neglecting the free tubulin concentration variations in the microtubule
dynamics process, the dynamic instability equations governing
the time evolution of
density functions $p^+(x,t)$ and $p^-(x,t)$
are:

\begin{equation}
\frac{\partial p^+(x,t)}{\partial t}=-\frac{\partial {(v^+ p^+(x,t))}}{\partial x}-f^-_+p^+(x,t)+ f^+_-p^-(x,t)
\end{equation}
\begin{equation}
\frac{\partial p^-(x,t)}{\partial t}=v^-\frac{\partial p^-(x,t)}{\partial x}+f^-_+p^+(x,t)- f^+_-p^-(x,t),
\end{equation}
where $p^{\pm}(x,t)$ is the density function of growing (shrinking) microtubules
 in the interval $[x, x+\delta x]$.
(Hence, $p^\pm/p$ is the probability density
function of growing (shrinking) in the interval $[x, x+\delta x]$,
where $p:=p^++p^-$ denotes total density.)

In prior work of Dogterom and Leibler, the main result was the prediction of a sharp transition between an unlimited growth, with the average speed $J>0$
and a steady-state or bounded growth, characterized by a MT length distribution with $J=0$. In the unbounded growth region the average length increase is $<L>=Jt$, where
\begin{equation}
J=\frac{v^+f^-_+-v^-f^+_-}{f^+_-+f^-_+}
\end{equation}
and the distribution approaches asymptotically to a Gaussian of width $s\sqrt{D_{eff}t}$, where
\begin{equation}
D_{eff}=\frac{f^+_-f^-_+}{(f^+_-+f^-_+)^3}(v^+ + v ^-)^2
\end{equation}
and

In the steady state the distribution of MT length is exponential with mean:

\begin{equation}
<L>=\frac{v^-v^+}{v^-f^+_--v^+f^-_+}
\end{equation}
\section{Model}\label{model}
\vspace{.3 in}
Following the Dogterom and Leibler model for growing/shrinking MT, we consider a generalized one-dimensional concentration-dependent model for microtubule growth and shrinkage:
\be\label{sys}
\ba
\frac{\partial p^+(x,t)}{\partial t}&=-\frac{\partial (v^+p^+(x,t))}{\partial x}-f^-_+p^+(x,t)+ f^+_-p^-(x,t)+d\frac{\partial^2 (p^+(x,t))}{\partial{x^2}}
\\
\frac{\partial p^-(x,t)}{\partial t}&=v^-\frac{\partial p^-(x,t)}{\partial x}+f^-_+p^+(x,t)- f^+_-p^-(x,t)+d\frac{\partial^2 (p^-(x,t))}{\partial{x^2}}
\\
\frac{\partial c(x,t)}{\partial t}&=-kc(x,t)+ v^- p^-(x,t)-v^+ p^+(x,t)+ D\frac{\partial^2(c(x,t))}{\partial{x^2}}.
\ea
\ee

\begin{figure}[h!]
\centering
\includegraphics[scale=.4]{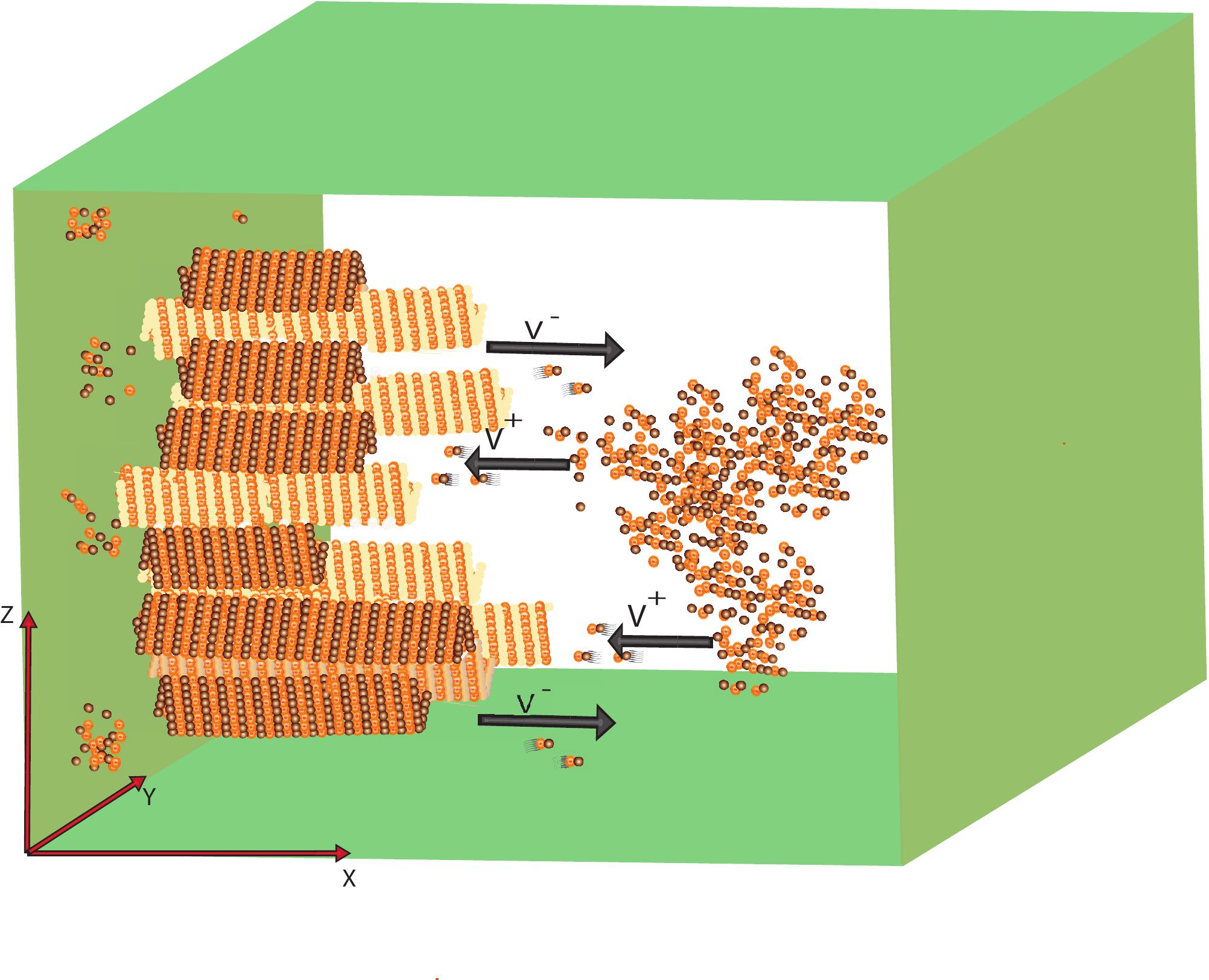}
\caption{Model geometry for infinitely rigid MTs originating from a planar surface}
\end{figure}

In these equations, $c(x,t)$ represents the concentration of the free
tubulin, $k$ is the hypothetic nucleation rate and
\be\label{cdependence}
f_-^+=\omega c(x,t) \; \hbox{\rm  and } \; v^+=u^+ c(x,t)
\ee
are the other concentration-dependent parameters.

The constant $k$ is assumed to be proportional to
the number of available nucleation sites per
volume; $kc$, measuring the rate of nucleation (with associated
loss of free subunits) is thus proportional to
the cross-section of a free subunit meeting a nucleation site,
the simplest possible assumption.
Likewise, \eqref{cdependence} are derived from
probabalistic cross-sections.
{\it Note that, properly, each of the coefficients $\omega$,
$u^+$, and $k$ should be assumed proportional to diffusion
constant $D$ measuring mean free path of subunits.}
As we hold $D$ fixed in this analysis, there is no harm in taking
them as constants however.

In this setup the nucleation points assumed hypothetically to be uniformly distributed throughout
$[0,+\infty)$. This model in the matrix form is as follows:

\be\label{eq:model}
U_t= A(U)U_x+ B(U)U+ CU_{xx}
\ee
on $x>0$, where
\be\label{U}
U(x,t)=
\bp
p^+(x,t)\\p^-(x,t)\\c(x,t)
\ep
\ee
\be\label{A}
A(U)=
\bp
-u^+c(x,t) & 0 &-u^+p^+(x,t)\\
0 & v^- & 0\\
0 & 0 & 0
\ep
\ee
\be\label{B}
B(U)=
\bp
-f^-_+ & \omega c(x,t) & 0\\
f^-_+ & -\omega c(x,t)& 0\\
- u^+c & v^- & -k
\ep
\ee
and
\be\label{C}
C=
\bp
d&0&0\\
0&d&0&\\
0&0&D&\\
\ep,
\ee
and $u^+$, $v^-$, $k$, $f^-_+$, $\omega$, $d$, $D$ are positive constants.
Here, $p^\pm\ge 0$ are
densities
of growth/decay of microtubules and
$c\ge 0$ is concentration of free tubulin.\\
Typical value for the parameters are
$$
D=0.5, \quad \omega=0.15, \quad \nu^-=0.05, \quad f_+^-=0.0005 ,\\ \quad d \simeq 0,
\quad
$$
with $u^+$, $k$ positive constants for which we don't know yet
the typical values.

Boundary conditions at $x=0$
are Dirichlet conditions on the
densities $p^\pm$,
\be\label{pbc}
\bp p^+\\p^-\ep|_{x=0}=
\bp \bar p^+_0\\\bar p^-_0\ep,
\ee
and on the concentration $c$ either a Dirichlet condition
\be\label{cbcdir}
c|_{x=0}=\bar c_0
\ee
or a homogeneous Neumann condition
\be\label{cbcneumann}
\partial_x c|_{x=0}=0.
\ee

\subsection{Goals}\label{s:goals}

We seek to study the existence and stability of {\it steady-state
solutions} of \eqref{eq:model}--\eqref{cbcneumann} of ``boundary-layer'' type,
that is, for which the solution $U$ approaches a constant state
$U_+$ as $x\to +\infty$.
Precisely, we seek asymptotically constant stationary solutions
\be\label{stat}
U(x,t)\equiv \bar U(x),
\quad \lim_{z\to +\infty} \bar U(z)=U_+,
\ee
of \eqref{eq:model},
i.e., asymptotically constant solutions of the steady-state ODE
\be\label{ode}
A(U)U_x+ B(U)U+ CU_{xx}=0
\ee
with boundary conditions \eqref{pbc} and \eqref{cbcdir} or \eqref{cbcneumann}
at $x=0$.
When such solutions exist, we seek to study their {\it stability}, that is,
whether a perturbation $\tilde U$ satisfying the same equation
\eqref{eq:model} that is close to $\bar U$ at initial time $t=0$
in some norm remains close to $\bar U$ for all $t>0$ in some
(possibly different) norm.
In applications, it is only stable steady states that are truly
steady in a practical sense;
unstable steady states, though steady in an idealized sense,
persist indefinitely only for a measure (hence probability) zero
set of initial data.

\section{Calculation of the endstates $U_+$}\label{equilibria}

Preparatory to the study of steady state solutions $\bar U$ of \eqref{stat},
we first investigate the possible limiting states
$U_+$ to which $\bar U$ may converge, seeking rest points
\be\label{eq:equil}
B(U_+)U_+=0
\ee
of \eqref{ode}.

The first two equations of \eqref{eq:equil} yield
$ p^-=\frac{f_+^-}{\omega c}p^+$, whereupon the third equation
yields
\be\label{psolve}
\Big( -\frac{u^+\omega c^2}{f_+^-}+ \nu^- \Big)p^-= ck.
\ee
Solving, we obtain the one-parameter family of solutions
\be\label{nonzeroeq}
\cE:=\{ p^-=\frac{kcf^-_+}{\nu^-f^-_+-u^+\omega c^2},\,
\, p^+=\frac{\omega k c^2}{\nu^-f^-_+-u^+\omega c^2}
,\, c \ge  0 \; \hbox{\rm arbitrary}\}.
\ee
From the requirement $0\le p^-<+\infty$, we
obtain the physicality condition
\be\label{genass}
\nu^-f^-_+-u^+\omega c^2 > 0.
\ee

\begin{rem}
\textup{
Equations \eqref{psolve}, \eqref{genass} simply
reflect the pair of balance laws $f^-_+p^+=f^+_-p^-$ and
$\nu^-p^- - \nu^+p^+=ck>0$ expressing conservation of
densities
and free subunit concentration, respectively,
in the absence of spatial variation.
Recall that $ck$ is the rate at which free subunit concentration is lost
to nucleation, or binding of subunits to nucleation sites.
}
\end{rem}

\subsection{Exponential decay to $U_+$}
We next investigate the rate of decay of solutions approaching
endstates $U_+$, specifically, whether or not approach at uniform
exponential rate.
As there is a continuous curve $\cE$ of viable endstates, this amounts
to verifying that the linearized ODE about $U_+$, written as a first-order
system, has a center subspace of dimension one, with no other neutral modes.

Linearizing \eqref{ode} about a rest point $U\equiv U_+$, we obtain
\be\label{lineq}
\tilde B(U_+) U + A(U_+)U_x+CU_{xx},
\ee
where
\be
\ba\label{tildeB}
\tilde B(U_+) &=B(U_+)+
\bp
0 & 0 & -\omega p^-_+\\
0 & 0 & \omega p^-_+\\
0 & 0 & -u^+p^+_+\\
\ep
=
\bp
-f^-_+ & \omega c_+& -\omega p^-_+\\
f^-_+ & -\omega c_+& \omega p^-_+\\
- u^+c_+ & \nu^- &  -k -u^+p^+_+\\
\ep.
\ea
\ee
It is convenient to find eigenvalues of \eqref{lineq} directly
in second-order form, substituting $U=e^{\mu x} V$ in \eqref{lineq},
to obtain
$$
(\mu A(U_+) + \tilde B(U_+) +\mu^2 C) V=0,
$$
or
\be\label{detcalc}
\ba
0&=\det( \mu A(U_+) + \tilde B(U_+) +\mu^2 C)\\
&=
\det \bp
d\mu^2-u^+c_+\mu-f_+^- & \omega c_+ & -\omega p^-_+ - u^+p^+_+ \mu\\
f_+^-& d\mu^2 +\nu^-\mu-\omega c_+  & \omega p^-_+\\
- u^+c_+ & \nu^- &   c\mu^2 -u^+p^+_+-k\\
\ep\\
&=
\mu
\det \bp
d\mu-u^+c_+ &   d\mu +\nu^- & -u^+p^+_+ \\
f_+^-& d\mu^2 +\nu^-\mu-\omega c_+  & \omega p^-_+\\
- u^+c_+ & \nu^- &   D\mu^2 -u^+p^+_+-k\\
\ep\\
&=\mu q(\mu), \\
\ea
\ee
\be\label{qdef}
\ba
q(\mu)&:=
(Dd^2)\mu^5+\Big(Dd\nu^--dDu^+c\Big)\mu^4\\
&-\Big(d^2u^+p^+_+ +d^2k+ f^+_-dD+ u^+c_+\nu^-D+d\omega c_+D\Big)\mu^3\\ &
+\Big(-f^+_-\nu^-D+ u^+c^2_+\omega D
+ u^+c_+dk- dv^-u^+p^+_+- d\nu^-k\Big)\mu^2 \\ &
+ \Big(d\omega c_+k -d\omega p^-_+\nu^- +u^+c_+\nu^-k + d\omega c_+u^+p^+_+ +f_+^-dk +f^-_+du^+p^+_+-u^+c_+d\omega p^-_+\Big)\mu
 \\ &
 +\Big(f_+^-\nu^-k-u^+c_+^2\omega k\Big),\\
\ea
\ee
where we have obtained the third equality in \eqref{detcalc}
by adding row two to row one and dividing out $\mu$ from the result.

So long as $q$ has no pure imaginary roots, then the center subspace is
dimension one and we obtain exponential decay.
We now check this, first for zero roots, then for nonzero imaginary roots.

\medskip
\subsection{ Checking zero roots}
Plugging in $\mu=0$,
we get $f_+^-\nu^--u^+c_+^2\omega=0$,
which is excluded by \eqref{genass}.
{\it Thus, we may conclude that there are no additional zero roots.}

\subsection{Checking pure imaginary roots}
Plugging in $\mu=i\xi$ in \eqref{detcalc} we get
\be\label{Imgroot}
\ba
q&(i\xi)=
\\&
(Dd^2)i\xi^5+\Big(Dd\nu^--dDu^+c\Big)\xi^4
+\Big(d^2u^+p^+_+ +d^2k+ f^+_-dD+ u^+c_+\nu^-D+d\omega c_+D\Big)i\xi^3\\ &
-\Big(-f^+_-\nu^-D+ uc^2_+\omega D
+ u^+c_+dk- dv^-u^+p^+_+- d\nu^-k\Big)\xi^2 \\ &
+ \Big(d\omega c_+k -d\omega p^-_+\nu^- +u^+c_+\nu^-k + d\omega c_+u^+p^+_+ +f_+^-dk +f^-_+du^+p^+_+-u^+c_+d\omega p^-_+\Big)i\xi
 \\ &
 +\Big(f_+^-\nu^-k-u^+c_+^2\omega k\Big)=\\&
  i\xi\Bigg(\Big(Dd^2\Big)\xi^4+\Big(d^2u^+p^+_+ +d^2k+ f^+_-dD+ u^+c_+\nu^-D+d\omega c_+D\Big)\xi^2 \\&
  +\Big(d\omega c_+k -d\omega p^-_+\nu^- +u^+c_+\nu^-k + d\omega c_+u^+p^+_+ +f_+^-dk +f^-_+du^+p^+_+-u^+c_+d\omega p^-_+\Big)\Bigg)+ \\&
  \Bigg( \Big(Dd\nu^--dDu^+c\Big)\xi^4-\Big(-f^+_-\nu^-D+ uc^2_+\omega D
+ u^+c_+dk- dv^-u^+p^+_+- d\nu^-k\Big)\xi^2\\
&+\Big(f_+^-\nu^-k-u^+c_+^2\omega k\Big)\Bigg).
\ea
\ee

We need only check nonexistence of pure imaginary roots
$\mu=i\xi$ with $\xi\ne 0$, $\xi$ real.
Noting that $q(i\xi)= i \xi q_1(\xi^2) + q_2(\xi^2)$,
with $q_1$ and $q_2$ quadratic, we find that such
$\xi$ must be of form $\xi=\pm \sqrt x$, where $x$ is a common
real positive root of $q_1$ and $q_2$.  The unique common root
is easily found by taking the resultant (taking a multiple eliminating
the quadratic term and solving the resulting linear equation,
then substituting this result into $q_1$ or $q_2$ to check whether
it vanishes- also checking whether it is positive).
This computation is carried out in Appendix \ref{s:resultant}.

Note that the resultant procedure described gives in the end an
analytic function of the parameters of the problem, which vanishes
if and only if $q_1$ and $q_2$ have a common root.  By properties
of analytic functions, it therefore vanishes either on a surface of
measure zero in parameter space, or else for arbitrary choice of parameters.
But, it is readily checked for specific parameters
that there is no common root (see Appendix \ref{s:resultant}).
Collecting information we may conclude the following general fact.

\begin{prop}\label{expdecay}
For generic choice of endstates $U_+\in \mathcal{E}$ (i.e., all
but a set of zero one-dimensional measure), all orbits converging
to $U_+$ do so at uniform exponential rate.
\end{prop}

\section{Dimension of the stable manifold at $U_+$}\label{dstab}

We next determine the dimension of the stable manifold at $U_+$,
that is, the number of eigenvalues $\mu$ of \eqref{detcalc} with
negative real part, or, equivalently, the number of roots of $q(\mu)=0$.
As a first step, computing the mod-two {\it stability index}
$$
\sgn\, q(0) q(+\infty)=
\sgn\,
(f_+^-\nu^-k-u^+c_+^2\omega k)(Dd^2)=+1
$$
using \eqref{genass}, we find that the number of unstable (i.e., positive
real part) roots is {\it even}, so that the number of stable
(negative real part) roots is {\it odd}, and thus is $1$, $3$, or $5$.

Next, using a standard technique in the stability theory,
we study the dispersion relation for the time-evolutionary
problem, from which we may conclude
by homotopy argument that the dimension is three.
Specifically, we look at the linearized time-evolutionary
equations about $U_+$,
\be\label{elineq}
\lambda U+ \tilde B(U_+) U + A(U_+)U_x+CU_{xx},
\ee
written as a first-order system, and look again at the number
of negative real part eigenvalues, or solutions $\mu$ of
the {\it indicial equation}
\be\label{indicial}
0=\det(\lambda U+ \mu A(U_+) + \tilde B(U_+) +\mu^2 C)
\ee
on the domain $\Lambda:=\{\lambda: \, \Re \lambda \ge 0\}$ of
interest for the later stability analysis.

Setting $\mu=i \xi$, $ \xi$ real, we obtain the {\it dispersion equation}
\be\label{disp}
0=\det(-\lambda U+ i \xi A(U_+) + \tilde B(U_+) -  \xi^2 C),
\ee
determining a family of three curves
\be\label{dispcurves}
\lambda_j( \xi)\in \sigma (i \xi A(U_+) + \tilde B(U_+) -  \xi^2 C),
\ee
$j=1,\dots, 3$,
where $\sigma(M)$ denotes the spectrum (eigenvalues) of a matrix $M$.
\\

%
%
%
%

We make the following standard assumption, verified in Section \ref{ver}
for $c_+$ small and checked numerically for the specific cases
treated in this paper.
\medskip

\begin{ass}\label{gooddisp}
The constant solution
$U\equiv U_+$ is linearly stable, i.e.,
\be\label{dispstab}
\Re \lambda_j( \xi)\le 0
\ee
for all $\xi\in \R$, with equality only for $\xi=0$.
\end{ass}

\begin{lem}\label{numbers}
Under Assumption \ref{gooddisp},
the numbers of stable and unstable roots of $\mu $
of \eqref{indicial} are independent of $\lambda$ for all
$\Re \lambda \ge 0$, $\lambda\ne 0$, and are both equal to three.
\end{lem}

\begin{proof}
The roots are continuous as functions of $\lambda$, so for the
first assertion we need only
check that no roots can cross the imaginary axis, i.e.,
there are no roots $\mu=i \xi$ with $ \xi$ real and $\Re \lambda\ge 0$
and $\lambda\ne 0$.
But, this possibility is ruled out by Assumption \ref{gooddisp}.
Thus, the numbers of stable and unstable roots are constant.
To determine their values, we may take $\lambda \to +\infty$
along the real axis,
and compute directly that in the limit there are three of each.
\end{proof}

\begin{cor}\label{c:three}
Under Assumption \ref{gooddisp},
the dimension of the stable manifold at states $U_+\in \mathcal{E}$
is generically three.
\end{cor}

\begin{proof}
By continuity of the roots of \eqref{indicial}
with respect to $\lambda$, together with the fact already established
that at $\lambda=0$
there is generically a single zero root,
we find, taking the limit as $\lambda \to 0$,
that there must be for $\lambda=0$ either two or three stable roots,
depending whether the single zero root is the limit of a stable or of
an unstable root as $\lambda \to 0$.
Since the number is odd, by our index computation, we conclude that
the dimension is generically three.
\end{proof}

\begin{rem}\label{incomingrmk}
\textup{
Note by our computations that we have obtained also the additional
information, which will be important for later stability analysis,
that the zero root at $\lambda=0$ perturbs as the real part of $\lambda$
is increased to the {\it unstable half plane}.
As we shall see later, this has the important consequence that the single
undamped mode governing ``total density'' $p:=p_++ p_-$
is convected inward toward the boundary $x=0$, and not away
toward $x=+\infty$, and therefore the nonlinear stability theory
may be treated by simpler weighted-norm methods \cite{He,Sat}
rather than the delicate pointwise methods of \cite{YZ,NZ}.
}
\end{rem}

%

\subsection{Verification for $c_+$ small}\label{ver}
We now verify Assumption \ref{gooddisp} for $c_+$ small, by explicit
computation of the case $c_+=0$.

\begin{prop}\label{smallc}
Assumption \ref{gooddisp} holds for
$U_+\in \mathcal{E}$ and $c_+$ sufficiently small.
\end{prop}

\begin{proof}
Plugging $U=e^{i\xi x}$  in \eqref{elineq}
and the fact that $\lambda \in \sigma(\tilde{B}+Ai\xi-C\xi^2)$,
we get
\be \label{stabcalc}
\ba
\\&
\tilde{B}+Ai\xi-C\xi^2=\\&
\bp
-f^-_+ & \omega c_+& -\omega p^-_+\\
f^-_+ & -\omega c_+& \omega p^-_+\\
- u^+c_+ & \nu^- &  -k -u^+p^+_+\\
\ep + i\xi\bp
-u^+c_+ & 0 &-u^+p^+_+\\
0 & \nu^- & 0\\
0 & 0 & 0
\ep -\xi^2\bp
d&0&0\\
0&d&0&\\
0&0&D&\\
\ep=
\\&
\bp
-f^-_+-i\xi u^+c_+-\xi^2d&\omega c_+& -\omega p^-_+-i\xi u^+p^+_+\\
f_+^-&-\omega c_++i\xi v^--\xi^2d&\omega p^-_+\\
-u^+c_+&\nu^-&-k-u^+p^+_+-\xi^2D
\ep
\ea
\ee
For case $\cE$, plugging
\be
p^+=\frac{\omega k c^2}{\nu^-f^-_+-u^+\omega c^2},
\quad
 p^-=\frac{kcf^-_+}{\nu^-f^-_+-u^+\omega c^2}
\ee
in \eqref{stabcalc}, we get
\be
\tilde{B}+Ai\xi-C\xi^2=
\bp
-f_+^--\xi u^+c_+i-\xi^2d & \omega c_+ & \frac{k\omega c_+(f_+^-+\xi u^+c_+i)}{-\nu^-f_+^-+u^+\omega c_+^2}\\
f_+^-&-\omega c_++\xi\nu^-i-\xi^2d&-\frac{k\omega c_+f_+^-}{-\nu^-f_+^-+u^+\omega c_+^2}\\
-u^+c_+&\nu^-&-\frac{k(-\nu^-f_+^--\xi D\nu^-f_+^-+\xi Du^+p^-_+c_+^2)}{-\nu^-f_+^-+u^+\omega c_+^2}
\ep
\ee
\medskip
When $c_+=0$, the eigenvalues are
\be\label{czero}
\lambda_1= -f_+^--\xi^2d,
\quad
\lambda_2=-\xi^2d+i\xi\nu^-,
\quad
\lambda_3=-k-\xi^2D,
\ee
which indeed have nonpositive real part for all real $\xi$.

In the general case, a straightforward computation shows that at $\xi=0$
there is a single zero eigenvalue (note that the lower lefthand
$2\times 2$ block of $\tilde B$ is nonsingular, by \eqref{genass})
of $\tilde B_+$, with
left and right eigenvectors $L=(1,1,0)$ and $R=(a,b,1)^T$,
where $a+b= \frac{f^-_+k}{\nu^- f^-_+-\omega c_+^2u^+}$ and
\be\label{pertform}
\alpha:=LA_+R= \frac{k}{a+b}
=\frac{\nu^- f^-_+-\omega c_+^2u^+}{f^-_+}>0,
\ee
whence, by standard matrix perturbation theory, the Taylor
expansion of this eigenvalue at $\xi=0$ is
$$
\lambda_2(\xi)= i \alpha \xi+\beta \xi^2+\dots+ \delta \xi^3,
$$
$\alpha>0$
verifying directly our observation of inward convection in the neutral mode.
(Note: the fact that $\alpha$ is real follows with no computation, simply
from \eqref{pertform}.)
Moreover, we may verify $\Re \lambda_1(0)<0$, $\Re \lambda_3(0)<0$
by computing the characteristic polynomial
$$
\det (\tilde B_+-\mu I)=\mu (\mu^2 + r \mu + s)
$$
and verifying the discriminant condition
$$
r^2-4s
=(u^+p^+_++k +\omega c_+ -f^-_+)^2
+4(f^-_+k  + u^+p^-_+ )>0.
$$

For $c_+$ sufficiently small,
$\Re \beta <0$, $\delta$ uniformly bounded,
by continuity of the Taylor coefficients in the limit.
This establishes that $\Re \lambda_2\le 0$ for $\xi$ near the origin,
with equality only at $\xi=0$.  On the other hand, $\Re \lambda_1$ and
$\Re \lambda_2$ by continuity are strictly negative on any bounded set
of $\xi$ and all three are strictly negative (almost by inspection)
in the limit as $|\xi|\to \infty$.
This establishes the assumption for $c_+$ sufficiently small.
\end{proof}

\subsection{Verification for general $c_+$}\label{bigc}
For general $c_+$, we verify Assumption \ref{gooddisp} numerically
in the course of our Evans function computations, via the following
elementary observation.

\begin{lem}\label{el}
For $D,d>0$, Assumption \ref{gooddisp} is equivalent to the
property that for $\lambda$ pure imaginary and $0<|\lambda|\le R$,
$R>0$ sufficiently large, the number of roots $\mu$ of indicial
equation \ref{indicial} with negative (resp. positive) real parts
is constant and equal to three.
\end{lem}

\begin{proof}
It is evident that for $|\xi|$ sufficiently large
$\Re \lambda_j(\xi) \le -\theta \xi^2$, some $\theta>0$.
Thus, if Assumption \ref{gooddisp} is violated, then
$\Re \lambda_j(\xi)$ is pure imaginary for some $0<|\xi|\le R$,
or, equivalently, there is a pure imaginary root $\mu=i\xi$
of the indicial equation \ref{indicial} for some pure imaginary
$\lambda$ with $0<|\lambda|\le R$.
But, this means that either the number of negative real part
roots or the number of positive real part roots must be less than
three, or else the total of all roots would be at least seven,
a contradiction.
\end{proof}

As described below, the numbers of stable and unstable roots $\mu$
of the indicial equation are checked at each step of our Evans function
computations, in particular, on an imaginary interval $[-i\hat R,i\hat R]$
with $\hat R$ larger than the value $R$ described above; see
Remark \ref{linkrmk}.


\section{Existence of steady-state solutions}\label{exist}
\subsection{General theory}\label{existgen}

Writing \eqref{ode} in first-order form as
\be
\bp U\\P \ep_x =
\bp
P\\
-C^{-1}(A(U)P-B(U)U
\ep
\ee
with $P:=U_x$,
we obtain a six-dimensional dynamical system.
The three boundary conditions at $x=0$ determine
a four-dimensional manifold of solutions (three free parameters
in the initial value, plus one dimension in the direction of spatial-evolution
for a particular initial value at $x=0$).

By Corollary \ref{c:three}, under Assumption \ref{gooddisp},
the stable manifold associated with
a rest point $U_+$ generically has dimension three.
Thus, in looking for stationary solutions satisfying the boundary conditions
at $x=0$ and converging as $x\to +\infty$ to a specified
endstate $U_+$, we seek
the intersection of a four-dimensional with a three-dimensional manifold,
{\it which should generically consist of a one-dimensional manifold,
or a finite union of distinct curves, corresponding
to distinct steady-state solutions.}

That is, by a dimensional count, connections seem to be possible
for all endstates $U_+\in \mathcal{E}$ that are stable as constant
solutions, as (numerically) all examples considered seem to be.


\begin{rem}\label{blt}
\textup{
The above suggests, by formal matched asymptotic expansion,
that behavior far from the boundary (distance $>> \max\{d,D\}$)
of solutions of \eqref{eq:model}
should be governed by the {\it inviscid} equation
$U_t=A(U)U_x + B(U)U$, with boundary condition $U\in \mathcal{E}$
at $x=0$, whether for Dirichlet or Neumann conditions imposed in
the full viscous equations, so long as the boundary layer (i.e.,
steady state solution connecting to $\mathcal{E}$) is stable.
See \cite{GMWZ} for related discussion in more general context.
}
\end{rem}

\subsection{Numerical determination of steady-state profiles}\label{prof}

To determine the profiles $\bar U$ of steady-state solutions,
we use a numerical boundary-value solver,
imposing boundary conditions at $x=0$ and projective boundary conditions
at $+\infty$, ensuring correct entry toward a specified endstate in $\cE$ along
its stable manifold.

Expanding \eqref{ode} we have,

\begin{equation}
\begin{split}
0&=-u^+cp_x^+-u^+p^+c_x-f_+^-p^++\omega cp^-+dp_{xx}^+,\\
0&= \nu^-p_x^-+f_+^-p^+-\omega cp^-+dp_{xx}^-,\\
0&=-u^+cp^++\nu^-p^--kc+Dc_{xx},
\end{split}
\end{equation}
 which yields the first order system
 \begin{equation}\label{yode}
 \begin{split}
 	\begin{pmatrix}y_1\\y_2\\y_3\\y_4\\y_5\\y_6 \end{pmatrix}_x=
F(Y):=
	\begin{pmatrix}y_2 \\ \frac{1}{d} (u^+y_5y_2+u^+y_1y_6+f_+^-y_1-\omega y_5y_3)\\
	y_4\\ \frac{1}{d}(\omega y_5y_3-\nu^-y_4-f^-_+y_1)\\
	y_6\\
	\frac{1}{D}(u^+y_5y_1-\nu^-y_3+ky_5)
	\end{pmatrix},
 \end{split}
 \end{equation}
with
\be\label{ydefs}
Y=(y_1,y_2,y_3, y_4,y_5,y_6)^T:=
(p^+, p^+_x, p^-, p^-_x,c,c_x)^T.
\ee

 The Jacobian at $+\infty$ is given by
 \begin{equation}\label{jacobian}
\frac{\partial F}{\partial Y}(Y_+)=
 	\begin{pmatrix}
		0&1&0&0&0&0\\
		\frac{f_+^-}{d}& \frac{u_+c_+}{d}& -\frac{\omega c_+}{d}&0&-\frac{\omega p_+^-}{d}&\frac{u^+p_+^+}{d}\\
		0&0&0&1&0&0\\
		\frac{-f_+^-}{d}&0&\frac{\omega c_+}{d}& -\frac{\nu^-}{d} &\frac{\omega p_+^-}{d}& 0\\
		0&0&0&0 &0&1\\
		\frac{u^+c_+}{D} &0&-\frac{\nu^-}{D} &0&
\frac{u^+p_+^++k}{D}&0
	\end{pmatrix},
 \end{equation}	
with $Y_+= (p^+_+, 0 , p^-_+, 0 ,c_+ ,0)^T$.

\subsection{Numerical method} \label{nummethod}
We solve \eqref{yode} as a two-point boundary-value problem
on $[0,M]$ with $M>0$ chosen sufficiently large.
At $x=0$, we impose the three boundary conditions
\be\label{ybc}
y_1(0)=p^+_0,
\quad
y_3(0)=p^-_0,
\; {\rm and }\;
y_5(0)=c_0
\; {\rm or }\;
y_6(0)=0.
\ee
At $x=M$, following the general approach of \cite{Be}, we
impose {\it projective boundary conditions}
\be\label{Mbc}
L\cdot (y_1(M)-p^+_+)=0
\quad
L\cdot (y_3(M)-p^-_+)=0
\quad
L\cdot (y_5(M)-c_+)=0,
\ee
where $L$ is any $3\times 6$ complex matrix that is full rank on the
center unstable subspace of the Jacobian $\frac{\partial F}{\partial Y}(Y_+)$.
We use the simple and well-conditioned
choice of $L$ consisting of rows spanning the unstable left
eigenspace of $\frac{\partial F}{\partial Y}(Y_+)$.

Following \cite{BHRZ,HLZ,CHNZ}, we implement this scheme
using \textsc{MATLAB}'s boundary-value solvers {\tt bvp4c},
{\tt bvp5c}, or {\tt bvp6c} \cite{HM}, which are adaptive Lobatto quadrature schemes and can be interchanged for our purposes.
The value of approximate spatial infinity $M$ are determined experimentally
by the requirement that the absolute error $|U(M)-U_+|$
be within a prescribed tolerance, say $TOL=10^{-3}$.
For rigorous error/convergence bounds for these algorithms,
see, e.g., \cite{Be}.

\subsection{Results} \label{numprofresults}

For each choice of endstates and parameters tested, we successfully
found steady-state profiles.  A typical example is shown in Figure 3
for endstate $U_+=(0,0,0)$ and parameters
$d=.1$, $D=1$, $\omega=1$, $\nu^-=1$, $f^-_+=1$, $u^+=1$, $\kappa=1$.
Other examples may be found in Figures 7--9,
superimposed on time-evolution plots for nearby perturbed solutions.

\begin{figure}[h!]
\centering
\includegraphics[scale=1]{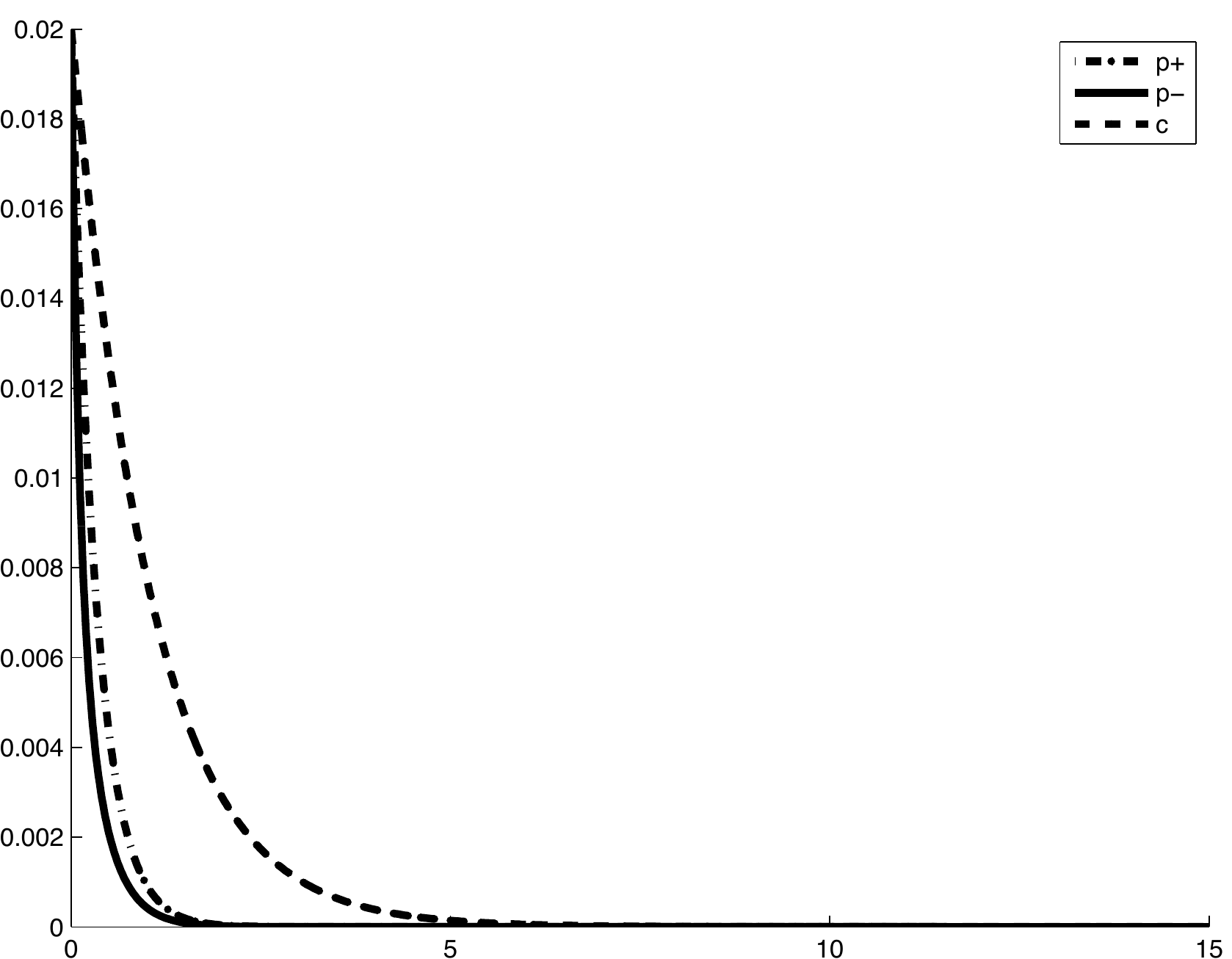}
\caption{Sample profile: $(p^+_+, p^-_+, c_+)=(0,0,0)$,
$d=.1$, $D=1$, $\omega=1$, $\nu^-=1$, $f^-_+=1$, $u^+=1$, $\kappa=1$.
 }
\end{figure}


\section{Abstract stability theory}\label{stab}
Under the assumption $D,d>0$,
$C$ is positive definite, so \eqref{eq:model}
is a system of semilinear second-order parabolic equations
of the type considered in \cite{He,Sat}.
Stability of steady-state solutions of this system may
be treated in straightforward fashion by the
{\it weighted norm method} of \cite{Sat}, as we now
briefly describe, {\it reducing the problem of determining
nonlinear stability to that of checking
spectral stability of the linearized operator about the wave.}
The simplicity of our treatment here
is to some extent an accident of the favorable structure
of this specific problem.  For (more complicated)
methods applying to general systems, see, e.g.,
\cite{HZ,Z2} and references therein.

\subsection{Abstract linearized equations}
Linearizing \eqref{eq:model} about a steady state $\bar U(x)$,
we obtain the {\it linearized equations}
\be\label{eq:lin}
U_t=LU:= A(x)U_x+ B(x)U+ CU_{xx},
\ee
$A(x):=A(\bar U(x))$, $B(x)V:=B(\bar U(x))V+ (dA(\bar U)V)\bar U_x$,
with boundary conditions
\be\label{linBC}
(p^+, p^-, c)=(0,0,0)
\; {\rm or }\;
(p^+,p^-,c_x)=(0,0,0).
\ee

\subsection{Abstract eigenvalue equations and spectrum}
The eigenvalue equations associated with \eqref{eq:lin} are
\be\label{eq:eig}
(L-\lambda)u=0
\ee
with boundary conditions \eqref{linBC}.
Let $W^{1,\infty}$ as usual denote the Banch space of bounded measurable
functions on $x\in [0,+\infty)$ possessing a bounded measurable weak
derivative, with norm $\|f\|_{W^{1,\infty}}:= \sup |f|+ \sup |f'|$.
Following \cite{He}, define the {\it spectrum} $\sigma(M)$ of an
operator $M$ with respect to a given Banach space $\calB$ to be the set
of $\lambda\in \CC$ for which $(\lambda-M)$ does not possess a bounded
inverse on $\calB$.
Define the {\it point spectrum} $\sigma_p(M)$
of $M$ to be the set of $\lambda \in \CC$
for which $(\lambda-M)U=0$ has a nonzero solution in $\calB$,
and the {\it essential spectrum} $\sigma_{ess}(M)$
of $M$ to be the set of $\lambda\in \CC$
that are in $\sigma(M)$ but not in $\sigma_p(M)$.

\begin{lem}\label{Helem}
Under Assumption \ref{gooddisp},
 $\Re \lambda\le 0$
for any $\lambda \in \sigma_{ess}(L)$, where
spectrum is defined with respect to $W^{1,\infty}$, with
equality at $\lambda=0$.
\end{lem}

\begin{proof}
By Lemma \ref{expdecay}, the coefficients of $L$ converge asymptotically
as $x\to +\infty$, whence, by a standard result of Henry \cite{He}
(adapted in straightforward fashion to the case of the half-line),
the essential spectrum of $L$ is either bounded to the left of the
rightmost envelope of the
dispersion curves $\lambda_j(\xi)$ of \eqref{disp}
or else contain all points to the right.
The latter possibility is easily eliminated by a standard energy estimate,
whence the result follows by \eqref{dispstab} and the fact that
$\lambda_2(0)=0$ (computed earlier).
\end{proof}

\subsection{The method of weighted norms}

Following \cite{Sat}, introduce now the weighted norm
\be\label{weight}
\|U\|_\eta:= \|e^{\eta x}U\|_{W^{1,\infty}},
\ee
$\eta\ge 0$, and the associated Banach space $\calB$
of functions with bounded $\|\cdot\|_\eta$ norm.
The following key lemma states that this
weighting has the effect of shifting the essential spectrum
of the linearized operator into the strictly stable
half-plane $\{\lambda: \Re \lambda<0\}$.

\begin{lem}\label{shift}
Under Assumption \ref{gooddisp}, for $\eta>0$ sufficiently small,
$\Re \lambda \le -\theta(\eta)$, $\theta>0$, for $\lambda\in \sigma_{ess}(L)$,
defined with respect to Banach space $\calB$.
\end{lem}

\begin{proof}
Making the change of variables $U\to e^{\eta x}U$
following \cite{Sat},
we convert $L$ to the operator
\begin{equation}\label{hatL}
\begin{aligned}
\hat L &:= C (\partial_x - \eta)^2 + A(\partial_x-\eta) + B
 =: \hat C \partial_x^2 + \hat A\partial_x + \hat B,
\end{aligned}
\end{equation}
where $\hat B=B+\eta A$, etc. (we need not compute $\hat A$ and $\hat C$).
Checking definitions, we find that
the spectrum of $L$ with respect to $\calB$ is the
spectrum of $\hat L$ with respect to $W^{1,\infty}$.

Applying the standard theory of \cite{He} as in the proof of Lemma \ref{Helem},
we find that $\sigma_{ess}(\hat L)$ with respect to $W^{1,\infty}$
is bounded on the left by the rightmost envelope of the dispersion
curves
$$
\hat \lambda_j(\xi)\in \sigma( \hat B_++ i\xi \hat A_+ -\xi^2C_+)
=\sigma \Big( B_++ (i\xi-\eta) A_+ +(i\xi-\eta)^2 C_+ \Big).
$$
Observing that $\lambda_j$ are continuous in $\eta$, we find by
Lemma \ref{Helem} that, for $1/R\le |\xi|\le R$,
$R>0$ arbitrary, $\Re \lambda_j<-\theta$, $\theta>0$ for $\eta>0$
sufficiently small.
Likewise, $\Re \lambda_j<-\theta$ for $|\xi|$ sufficiently large
and $\eta$ sufficiently small, by domination of term $C\xi^2$.

Referring now to the computations in the proof of Proposition \ref{smallc},
and following the notation therein,
we recall that $\Re \lambda_1$ and $\Re \lambda_3$ are strictly
negative for all $\xi$, so that the above continuity argument
in fact gives $\Re \hat \lambda_1, \, \Re \hat \lambda_2<-\theta$
for $\eta>0$ sufficiently small.

Hence, it remains only to check the behavior of $\hat \lambda_2$
for $|\xi|\le 1/R$ arbitrarily small.
Recalling the Taylor expansion \eqref{pertform} of $\lambda_2(\xi)$,
and substituting $\hat \lambda_2(\xi)=\lambda_2(\xi- \eta/i)$,
we find that
$$
\Re \hat \lambda_2(0)=
\Re i \alpha (\xi- \eta/i) + h(\xi- \eta/i )^{2j}) + O(|\xi-\eta/i|^{2j+1})
= -\alpha \eta + \Re h \xi^{2j} + o(|\xi|^{2j}+ \eta)
,
$$
where $ \alpha:=LAR >0 $ and $h$ is the first term with nonzero real
part in the Taylor expansion of $\lambda_2$ after $\alpha$, by
Assumption \ref{gooddisp} necessarily
an even-order term with $\Re h<0$.
Thus, $\Re \hat \lambda_2\le -\theta<0$ for some $\theta<0$,
for $|\eta,\xi|$ sufficiently small, completing the proof.
\end{proof}

\subsection{Basic nonlinear stability theorem}\label{s:stabthm}

We make a final assumption on stability of point spectrum of $L$.

\begin{ass}\label{stabcond}
$\Re \sigma_{p}(L) <0$.
\end{ass}

\begin{theo}\label{stabthm}
Under Assumptions \ref{gooddisp} and \ref{stabcond},
for $\eta>0$ sufficiently small,
let $\tilde U(x,t)$ be a solution of \eqref{eq:model}
with initial data $\tilde U_0$ such that
$\|e^{\eta x} (\tilde U_0-\bar U)\|_{W^{1,\infty}}$ is sufficiently small.
Then,
$\|e^{\eta x} (\tilde U-\bar U)\|_{W^{1,\infty}}(t)$
decays to zero as $t\to \infty$, at rate
$$
\|e^{\eta x} (\tilde U-\bar U)\|_{W^{1,\infty}}(t) \le Ce^{-\theta t},
$$
where
$\theta=\theta(\eta) >0$.
(Here, $\theta(\eta)\sim \alpha \eta $ as $\eta\to 0$, where $\alpha$
is as in \eqref{pertform}.)
\end{theo}

\begin{proof}
By Assumptions \ref{gooddisp} and Lemma \ref{shift},
we have $\Re \sigma_{ess}(L)\le -\theta<0$ with
respect to space $\calB$.
On the other hand, eigenvalues of $L$ with respect to
$\calB$ are necessarily eigenvalues with respect to $W^{1,\infty}$
as well, hence, by Assumption \ref{stabcond},
$\Re \sigma_{p}(L) \le -\theta<0$ with respect to
$\calB$ as well, and so
\be\label{gap}
\Re \sigma(L)\le -\theta<0
\ee
with respect to the space $\calB$.

By the spectral gap \eqref{gap} and the fact that $L$ as as second-order
elliptic operator is sectorial,
$L$ by standard analytic semigroup theory \cite{He}
is {\it linearly exponentially stable} with respect to $\calB$,
i.e.,
$$
\|e^{Lt}f\|_\eta \le e^{-\theta_1 t}\|f\|_{\eta}
$$
Turning to the nonlinear equations, and noting that $e^{\eta x}$
is an exponentially {\it growing} weight, we find as in
\cite{Sat} that nonlinear exponential stability follows as well,
by a standard nonlinear iteration/variation of constants argument.
We omit these straightforward details.
\end{proof}

Theorem \ref{stabthm} asserts that, similarly as for finite-dimensional
ODE, exponential stability of steady-state solutions of the PDE
\eqref{eq:model} reduces to spectral properties of the linearized
operator $L$: specifically, verification of Assumptions
\ref{gooddisp} and \ref{stabcond}.
This abstract result
reduces the study of PDE stability to the study of the
eigenvalue equation for $L$, an ODE.


\begin{rem}\label{gen}
\textup{
The success of the weighted norm method relies on the property
that the undamped mode $p=p^++p^-$ is convected inward toward
the boundary $x=0$.  Intuitively, we may think of
$p$ to lowest order as $\tilde p(x+\alpha t)$ for a fixed profile
$\tilde p$, so that
$$
\sup e^{\eta x}p=e^{-\eta \alpha t}
\sup e^{\eta(x+\alpha t)}\tilde p(x+\alpha t)
$$
decays time-exponentially for $\eta>0$.
More general situations may be treated by pointwise methods
as in \cite{HZ,Z2}, to obtain time-algebraic rather than
exponential rates of decay.
}
\end{rem}

\section{Spectral stability analysis}\label{specstab}
Theorem \ref{stabthm} reduces the study of nonlinear stability
to verification of spectral stability assumptions \ref{gooddisp} and
\ref{stabcond}, the first a somewhat standard linear algebraic problem
and the second a type of ordinary differential boundary-value problem
arising frequently in the study of eigenvalues of differential operators.
We treat these numerically using {\it numerical Evans function techniques}
developed in \cite{Br,BDG,HuZ}.

\subsection{Linearized eigenvalue equations}

Expanding \eqref{eq:lin} coordinate-wise, we obtain
the {\it linearized equations}
\be
\ba
\partial_t p^+&=
-u^+ \bar c  p^+_x -u^+ \bar c_x  p^+
-u^+ c  \bar p^+_x -u^+ c_x  \bar p^+
-f^-_+ p^+ +\omega \bar c p^- + \omega c \bar p^- + d p^+_{xx}
\\
\partial_t p^-&=
\nu^-  p^-_x
+f^-_+ p^+ -\omega \bar c p^- - \omega c \bar p^- + d p^-_{xx}
\\
\partial_t c&=-kc+ \nu^- p^-
- u^+ \bar c p^+ - u^+ c \bar p^+ + Dc_{xx},
\\
\ea
\ee
with boundary conditions $(p^+,p^-,c)=(0,0,0)$ (Dirichlet)
or $(p^+,p^-,c_x)=(0,0,0)$ (Neumann).

Seeking normal modes $U(x,t)=e^{\lambda t}u(x)$,
we obtain the {\it linearized eigenvalue equations}
\be\label{linevalue}
\ba
\lambda p^+&=
-u^+ \bar c  p^+_x -u^+ \bar c_x  p^+
-u^+ c  \bar p^+_x -u^+ c_x  \bar p^+
-f^-_+ p^+ +\omega \bar c p^- + \omega c \bar p^- + d p^+_{xx}
\\
\lambda p^-&=
\nu^-  p^-_x
+f^-_+ p^+ -\omega \bar c p^- - \omega c \bar p^- + d p^-_{xx}
\\
\lambda c&=-kc+ \nu^- p^-
- u^+ \bar c p^+ - u^+ c \bar p^+ + Dc_{xx}.
\\
\ea
\ee
with boundary conditions
\be\label{linBC2}
(p^+,p^-,c)=(0,0,0)
\; {\rm or }\;
(p^+,p^-,c_x)=(0,0,0).
\ee
We seek growing or neutral modes, i.e., bounded solutions of
\eqref{linevalue}--\eqref{linBC2} with $\Re \lambda \ge 0$.

\subsection{The Evans function}\label{evans}
Written as a first-order system, the eigenvalue equations
\eqref{linevalue} appear as
\be\label{evansode}
W_x=\mA(x,\lambda) W,
\ee
where
\be\label{wvar}
W=(w_1, w_2,w_3,w_4, w_5, w_6)^T:=
(p^+,p^+_x, p^-,p^-_x,c,c_x)^T
\ee
and
\be\label{mA}
\mA:=
\begin{pmatrix}
0 & 1 & 0 & 0 & 0 & 0\\
\frac{\lambda + u^+ \bar c_x + f^-_+}{d}& \frac{u^+ \bar c}{d} &
-\frac{ \omega \bar c}{d} & 0 &
\frac{u^+\bar p^+_x - \omega \bar p^-}{d} & \frac{u^+ \bar p^+}{d}\\
0 & 0 & 0 & 1 & 0 & 0\\
-\frac{f^-_+}{d} & 0 &
\frac{\lambda +\omega \bar c}{d}& -\frac{\nu^-}{d} &
\frac{\omega \bar p^-}{d}& 0\\
0 & 0 & 0 & 0 & 0 & 1\\
\frac{ u^+\bar c}{D} & 0 &
- \frac{ \nu^-}{D} & 0 &
\frac{ \lambda + u^+ \bar p^+ + k}{D}  & 0 \\
\end{pmatrix},
\ee
$x\in [0,+\infty)$,
with boundary conditions \eqref{linBC2} imposed at $x=0$.

The limiting coefficient matrix as $x\to +\infty$ is
\be\label{mAplus}
\mA_+(\lambda):=
\begin{pmatrix}
0 & 1 & 0 & 0 & 0 & 0\\
\frac{\lambda +  f^-_+}{d}& \frac{u^+ c_+}{d} &
-\frac{ \omega c_+}{d} & 0 &
-\frac{ \omega p^-_+}{d} & \frac{u^+ p^+_+}{d}\\
0 & 0 & 0 & 1 & 0 & 0\\
-\frac{f^-_+}{d} & 0 &
\frac{\lambda +\omega c_+}{d}& -\frac{\nu^-}{d} &
\frac{\omega p^-_+}{d}& 0\\
0 & 0 & 0 & 0 & 0 & 1\\
\frac{ u^+c_+}{D} & 0 &
- \frac{ \nu^-}{D} & 0 &
\frac{ \lambda + u^+ p^+_+ + k}{D}  & 0 \\
\end{pmatrix},
\ee
which for $\lambda=0$ reduces to the Jacobian
$\frac{\partial F}{\partial Y}(Y_+)$ computed in \eqref{jacobian}.
By exponential convergence of $\bar U(x)$ to $U_+$ as $x\to \infty$,
we have exponential convergence also of $\mA(x,\lambda)$ to $\mA_+$.

By Lemma \ref{numbers} and Corollary \ref{c:three},
under Assumption \ref{gooddisp},
the stable subspace of $\mA$ has dimension
three for $\Re \lambda\ge 0$, whence, by
a standard lemma\footnote{
The ``gap'' or ``conjugation'' lemma; see, e.g., \cite{GZ,Z1,GMWZ}.}
 using the exponential convergence of $\mA$
to $\mA_+$, there exist choices of three independent solutions
$W_j^+$, $j=4,5, 6$ of \eqref{evansode} spanning the
set of solutions that are decaying as $x\to +\infty$ that
are analytic as functions from $x\to W^{1,\infty}[0,+\infty)$.
Likewise, one may easily construct analytic choices of three independent
solutions $W_j^0$, $j=1,2,3$ of \eqref{evansode} satisfying the
boundary conditions \eqref{linBC} at $x=0$.
Evidently, $\lambda$ is an eigenvalue, i.e., there exists a solution
of \eqref{evansode} that is decaying as $x\to +\infty$ and satisfies
boundary-conditions \eqref{linBC} at $x=0$, if and only if there
is a linear dependency among
$W_1^0, W_2^0, W_3^0, W_4^+, W_5^+, W_6^+$.

Defining the {\it Evans function}
\be\label{evansfn}
E(\lambda):=
\det ( W_1^0, W_2^0, W_3^0, W_4^+, W_5^+, W_6^+)|_{x=0},
\ee
following \cite{AGJ,GZ,Z1}, we thus have that the eigenvalues
of $L$ with respect to the weighted space $\calB$
on $\Re \lambda\ge 0$ correspond to the zeros of $E$.
This can be efficiently computed numerically using the STABLAB package
developed by J. Humpherys, based on algorithms of \cite{BrZ,HuZ}.
We refer the reader to \cite{HuZ,HLyZ} for a discussion
of theory and numerical protocol, which is by now standard.\footnote{
See, e.g., \cite{BHRZ,HLZ,CHNZ,HLyZ,BHZ,BLZ,BLeZ},}

%

\subsection{High-frequency bound}

Define $\delta:=\min\{d,D\}$, $\alpha= \max |A(x)|$, $\beta=\max |B(x)|$,
$A,B,C$ the coefficients of $L$, where
$|M|$ denotes the Euclidean matrix operator norm of a matrix $M$,
i.e., the square root of the largest eigenvalue of $M^TM$.
\begin{lem}\label{hflem}
There exist no spectra of $L$ with respect to $\calB$ for
$\Re \lambda +|\Im \lambda| > \alpha^2/\delta + \beta$
and $\eta>0$ sufficiently small.
\end{lem}

\begin{proof}
First note that eigenvalues of $L$ with respect to $\calB$ are eigenvalues
with respect to $L^2$ as well.
Taking the real part of the $L^2$ complex inner product of $U$
with eigenvalue equation
\be\label{eq:eig2}
\lambda U= AU_x+ BU+ CU_{xx},
\ee
we obtain after an integration by parts of term $\int U^* CU dx$
the inequality
$$
\Re \lambda \|U\|_{L^2}^2
\le \alpha \|U\|_{L^2}^2 + \beta
\|U\|_{L^2} \|U_x\|_{L^2} -\delta \|U_x\|_{L^2}^2.
$$
Likewise, taking the imaginary part yields
$ |\Im \lambda | \|U\|_{L^2}^2 \le \alpha \|U\|_{L^2}^2 .  $
Summing, and using Young's inquality,
$2\alpha ab\le (\alpha^2/\delta)a^2 + \delta b^2$, yields
$$
(\Re \lambda +|\Im \lambda|) \|U\|_{L^2}^2
\le (\alpha^2/\delta +\beta) \|U\|_{L^2}^2 ,
$$
yielding
$\Re \lambda +|\Im \lambda|\le \alpha^2/\delta +\beta$
whenever $\|U\|_{L^2}\ne 0$.

This establishes nonexistence of point spectra for
$\Re \lambda +|\Im \lambda|> \alpha^2/\delta +\beta$.
Nonexistence of essential spectra follows by the same
argument applied to the linearized eigenvalue equations about the
constant solution $U\equiv U_+$ together with the observation
\cite{He} that essential spectrum of $L$ with respect to $L^2$ is
bounded on the left by the rightmost envelope of the spectrum
of the linearized operator about $U\equiv U_+$,
given by the solution set of the dispersion relation \eqref{disp},
and the fact that nonexistence of spectrum with respect to $\calB$
is implied by nonexistence of $L^2$ spectrum provided $\eta>0$ is
taken sufficiently small.
(See the argument of Lemma \ref{shift}.)
\end{proof}

\begin{rem}\label{linkrmk}
\textup{
The above energy estimate implies in passing
that $\hat R:= \alpha^2/\delta +\beta >R$, where $R$ is as in Lemma \ref{el}.
To verify Assumption \ref{gooddisp}, therefore, it is
sufficient to check that the stable subspace of $\mA_+$ has dimension three
for $\lambda$ on imaginary interval $[-i\hat R,i\hat R]$.
}
\end{rem}

\subsection{Numerical stability computations}

Following a standard strategy introduced by Evans and Feroe \cite{EF},
we may now check by a single {\it numerical winding number computation}
both of Assumptions \ref{gooddisp} and \ref{stabcond}.
Specifically, we compute the image under $E$ of a semicircle $S$ of radius
$\hat R:=\alpha^2/\delta +\beta$ (see Lemma \ref{hflem} for definitions)
in the unstable complex half-plane $\Re \lambda \ge 0$, with
diameter $[-i\hat R, i\hat R]$ lying along the imaginary axis.

As the code for approximating $E$ in the course of the computation
checks that the stable subspace of $\mA_+$ has dimension three,
we obtain by Remark \ref{linkrmk} a check
of Assumption \ref{gooddisp} in the course of computing $E$
on the diameter $\lambda\in [-i\hat R, i\hat R]$.
Assuming that Assumption \ref{gooddisp} is valid, we then obtain
by Lemma \ref{hflem} that there are no eigenvalues of $E$ outside
the semicircle.

Moreover, by the properties of the Evans function
described in Section \ref{evans} (recall: also dependent on Assumption
\ref{gooddisp}), eigenvalues of $L$ within the semicircle $S$ are exactly
the zeros of the analytic function $E$.
By the {\it Principle of the Argument}, the number of zeros within $S$
is equal to the winding number of $E(S)$, i.e., the number of times
$E(S)$ circles the origin in counterclockwise direction as $\lambda$
traverses $S$ in counterclockwise direction.
Thus, {\it Assumption \ref{stabcond} corresponds to winding-number zero}.

\subsection{Results}
Figure 4 displays the image of the semicircle $S$
under the (analytic) Evans function $E$ for the typical profile $\bar U$
displayed in Figure 3.
We see clearly
that the winding number is zero, indicating {\it spectral stability},
hence, by the results of Section \ref{stab}, $\bar U$ is
nonlinearly stable.
Here, $d=.1$, $D=1$, $\omega=1$, $\nu^-=1$, $f^-_+=1$, $u^+=1$, $\kappa=1$,
so that $\alpha=1$, $\beta=1.4343$, and $\delta=.1$,
and $\hat R=\alpha^2/\delta+\beta \approx 11.43$; the radius of
$S$ is taken as $\tilde R=12>\hat R$, following Lemma \ref{hflem}.
Though we did not carry out a systematic study over all parameter
values and profiles, {\it for all profiles checked, we obtained
results of zero winding number consistent with stability}.

\begin{rem}
\textup{
Stability of ``small-amplitude'' steady states
$\bar U$ with $\max |\bar U-U_+|$ sufficiently small,
could in principle be determined by a study of stability of constant solutions,
as described in a more general setting in \cite{GMWZ}.
However, this does not seem particularly interesting for applications,
and so we do not attempt to carry out such an analysis here.
}
\end{rem}

\begin{figure}[h!]
\centering
\includegraphics[scale=1]{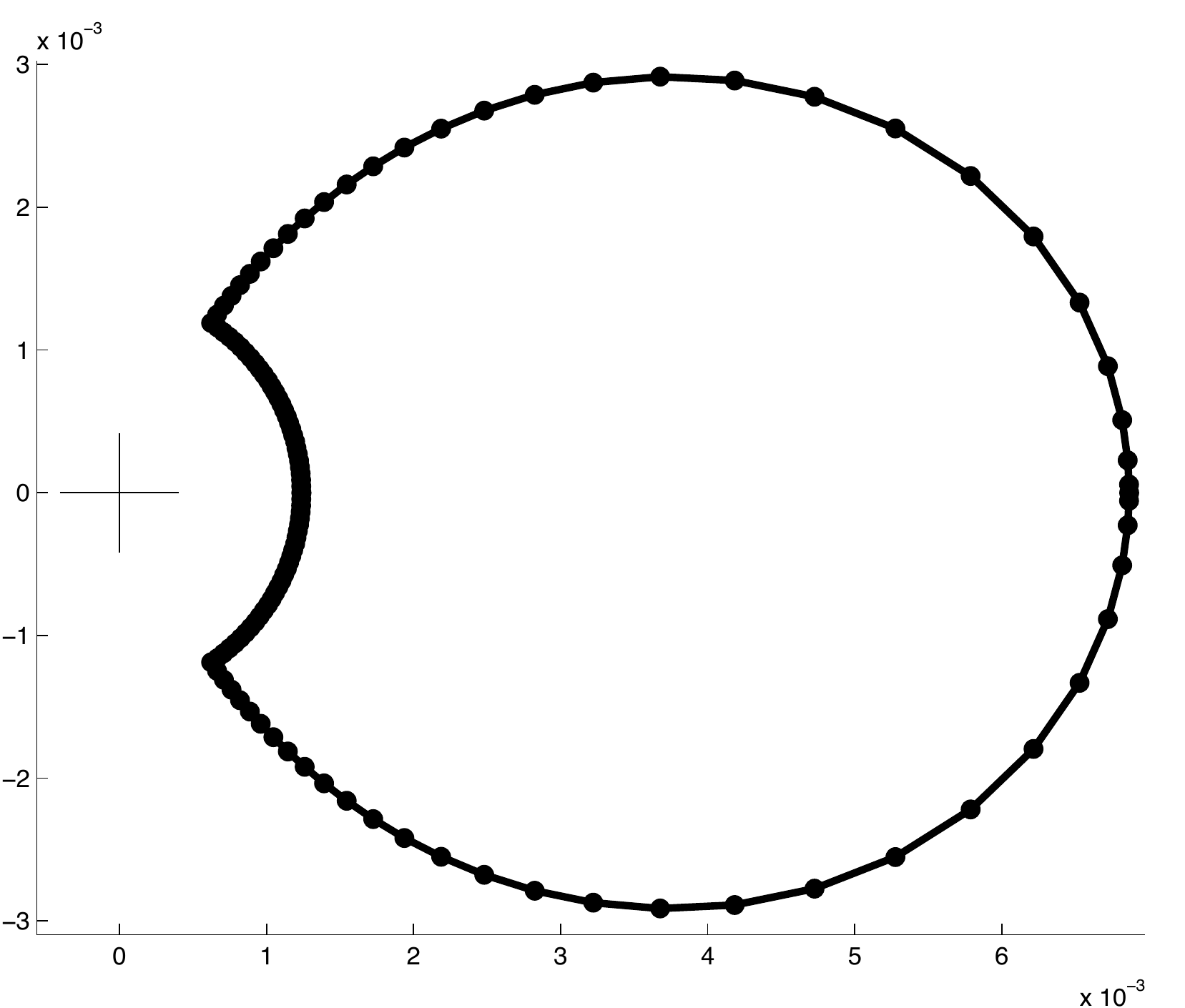}
\caption{Typical winding number computation,
corresponding the the profile of figure 3.
the displayed contour is the image under $D$ of
a semicircle of radius $12$, whereas
$\hat R\approx 11.43$.
 }
\end{figure}

\section{Behavior for general initial data}\label{behavior}

The results of Sections \ref{stab} and \ref{specstab} indicate
time-exponential stability of steady state solutions $\bar U$ with respect
to spatially-exponentially decaying initial perturbations.
To put this another way, solutions with initial data converging
exponentially as $x\to \infty$ to $U_+\in \mathcal{E}$ and
sufficiently close to the profile $\bar U$, will converge as $t\to \infty$
to $\bar U$.

However, as we explore in this section, the actual stability properties
appear to be considerably stronger.
Namely, initial data converging as $x\to\infty$
to {\it any } limiting value $U_\infty$ appears to settle down
rapidly to a steady state, determined by the total
density $p:=p^++p^-$ of the limiting state $U_\infty$.
(Indeed, we suspect that only $p$ need converge to a limiting
value in order to determine the limiting steady state.)

The weighted norm methods used above do not appear sufficiently
fine to establish this property, or at least we were not able
to carry out such an analysis.  An interesting open problem would
be to attempt to establish this result using the more detailed
pointwise Green function methods of \cite{HZ,Z1}.

%

\subsection{Numerical time-evolution approximations}\label{numtimeev}

In this section, we describe the results of a numerical time-evolution
study, that is, the evolution in time of given initial data under
equation \eqref{sys}.
This was carried out via MATLABs Finite Element Method (FEM) routines, with
Neumann boundary conditions at a boundary $x=50$ far from the computational
domain of interest.  For reference, see the helpful notes \cite{H}.


In Fig. 5, we have displayed graphs of variables $p^+$, $p^-$, and
$c$ taken at successive time intervals, with boundary data
$p^+_0$, $p^-_0$, and $c_0$ at $x=0$ fixed at the same values
$.2$, $.2$, $.2$ used in the computations used to obtain the
steady-state solution from Fig. 3, and step initial data consisting
of $.2$, $.2$, $.2$ for $x\in [0,3]$ and $0$, $0$, $0$ elsewhere.
We use the same model parameters as for Fig. 3 as well.

We see that the solution rapidly smooths and settles to an apparently
steady state represented by the lower envelope of the graphs,
with the darker bands near the lower envelopes indicating
convergence of each variable toward these curves.
In Fig. 6, we have superimposed onto Fig. 5 the steady-state profiles
that were displayed in Fig. 3, computed as solutions of a
two-point ordinary differential boundary problem,
showing a nearly exact correspondence
between the apparent limit of the time-evolution PDE solution
and the theoretical steady state ODE solution computed by completely
different methods.

This is strong evidence for the accuracy of both computations,
and also for stability of the steady-state solution as computed by
yet a third different method in our numerical Evans function computations.

\begin{figure}[h!]\label{evolfig}
\centering
\includegraphics[scale=.3]{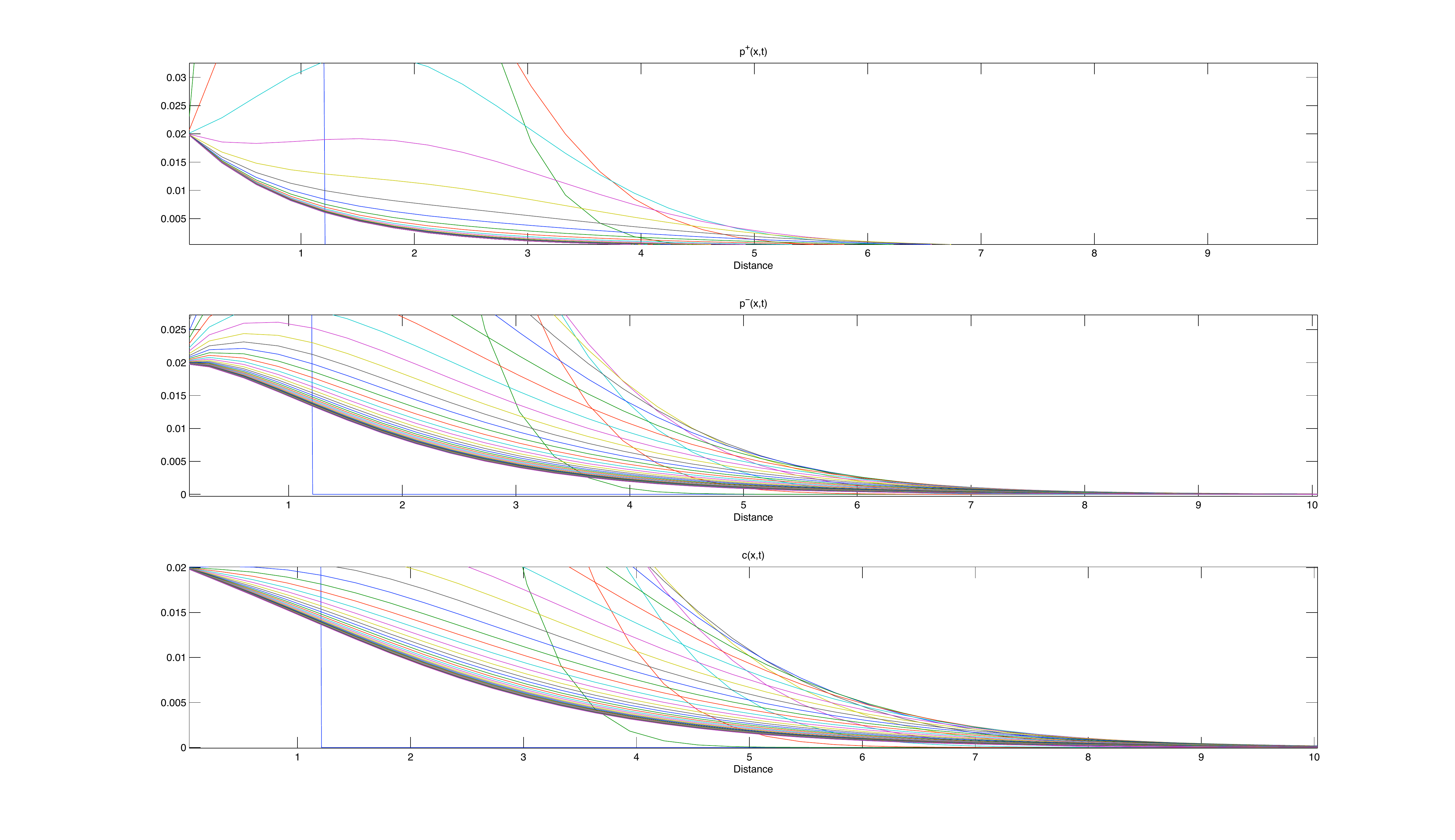}
\caption{Time evolution of typical data.
 }
\end{figure}

\begin{figure}[h!]
\centering
\includegraphics[scale=.35]{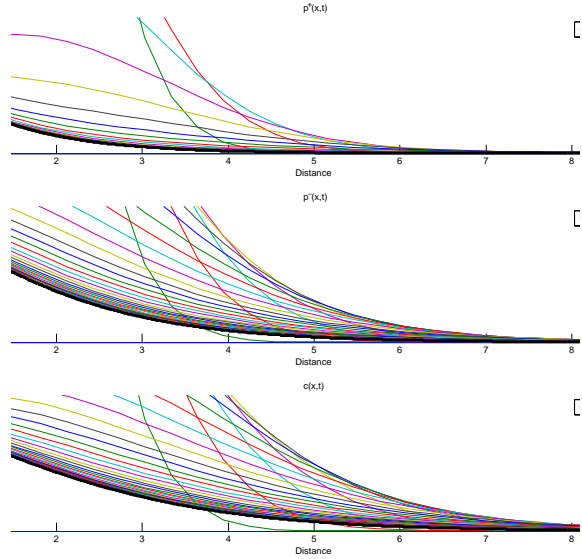}
\caption{Comparison with steady state profile.
 }
\end{figure}

\subsection{Selection of limiting profile}\label{selection}
The intuitive explanation for the behavior described in the previous
subsection is that only the special mode $p:=p^++p^-$ is ``undamped'',
with other modes decaying exponentially toward their equilibrium values;
see the discussion in Section \ref{equilibria}
of stability of constant solutions.
Thus, even if the limiting value $U_\infty$ of the initial data
is not on the equilibrium curve $\calE$,
it will rapidly adjust dynamically toward a value $U_+\inf \calE$
with the same total density $p$.

This is why the data from Fig. 5 converges toward a zero endstate profile.
In Fig. 7, we show the results of a similar experiment with the same
boundary data and the same
initial data for $p^\pm$, but $c$ initially constant and equal to $.02$;
as predicted by theoretical considerations, we see convergence to an
identical zero endstate profile.
In Fig. 8, we display the results for initial data $c\equiv .02$
and $p$ initially zero as $x\to + \infty$,
but $p^\pm$ nonzero as $x\to + \infty$
(to do this, we assign $p^-$ a nonphysical negative value, but
mathematically this still makes sense).

Further confirmation is given by Fig. 9, in which we display behavior
for initial data converging to nonzero value of $p$; again, the agreement
with prediction is exact.

\begin{figure}[h!]
\centering
\includegraphics[scale=.35]{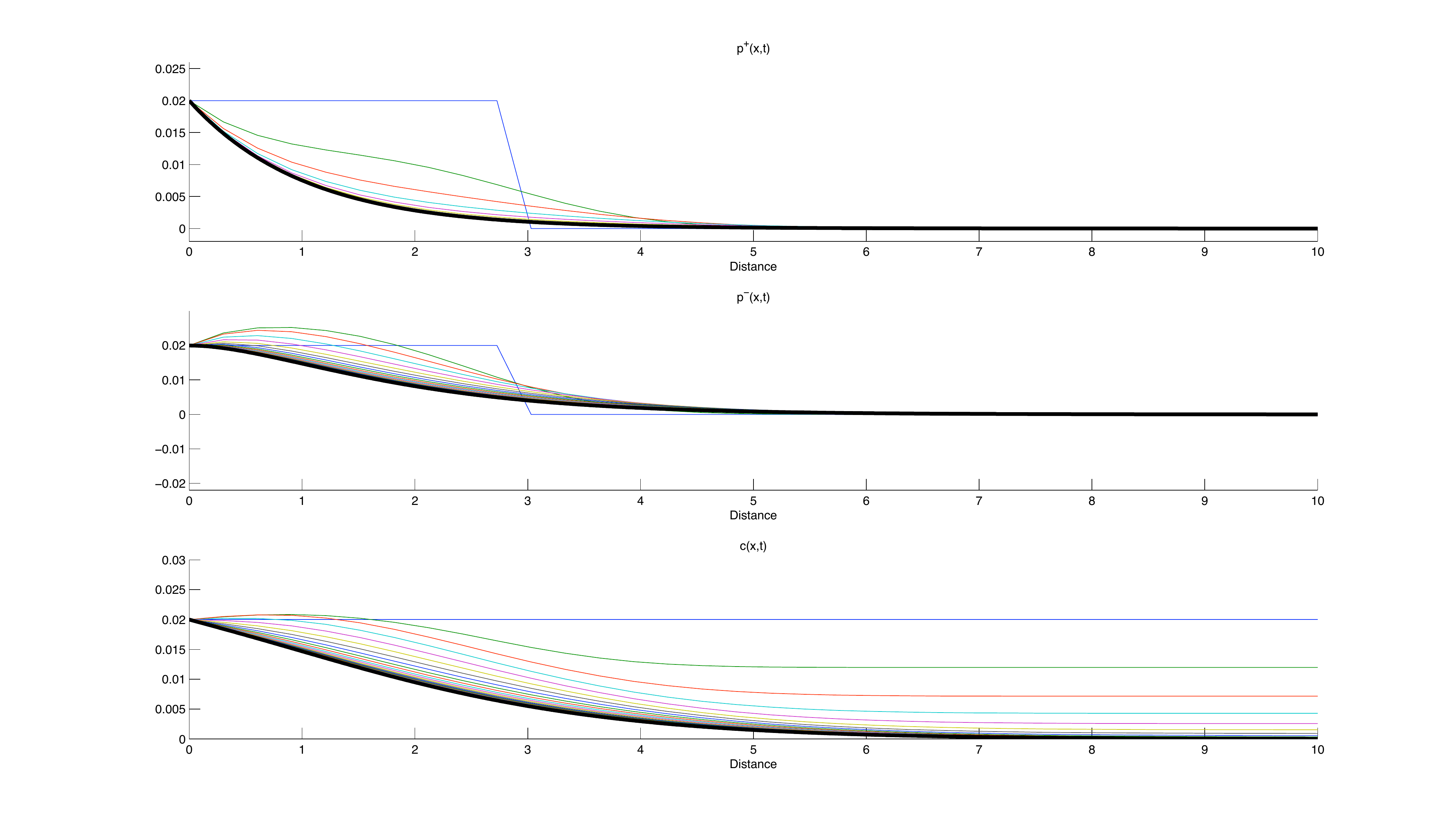}
\caption{Selection of endstate (i) $c$ initially nonzero.
 }
\end{figure}

\begin{figure}[h!]
\centering
\includegraphics[scale=.35]{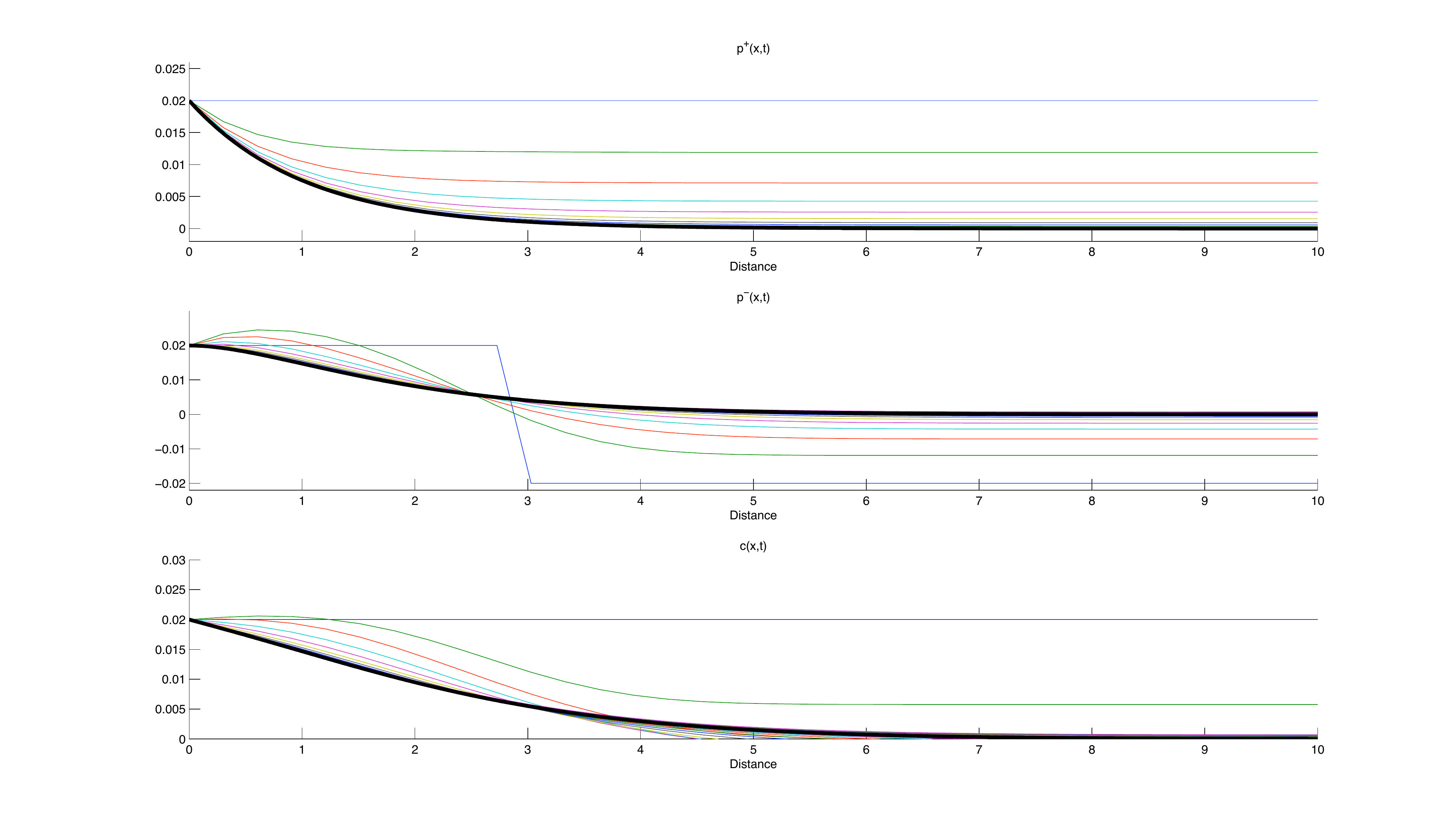}
\caption{Selection of endstate (ii) $p^\pm$ initially nonzero.
 }
\end{figure}

\begin{figure}[h!]
\centering
\includegraphics[scale=.35]{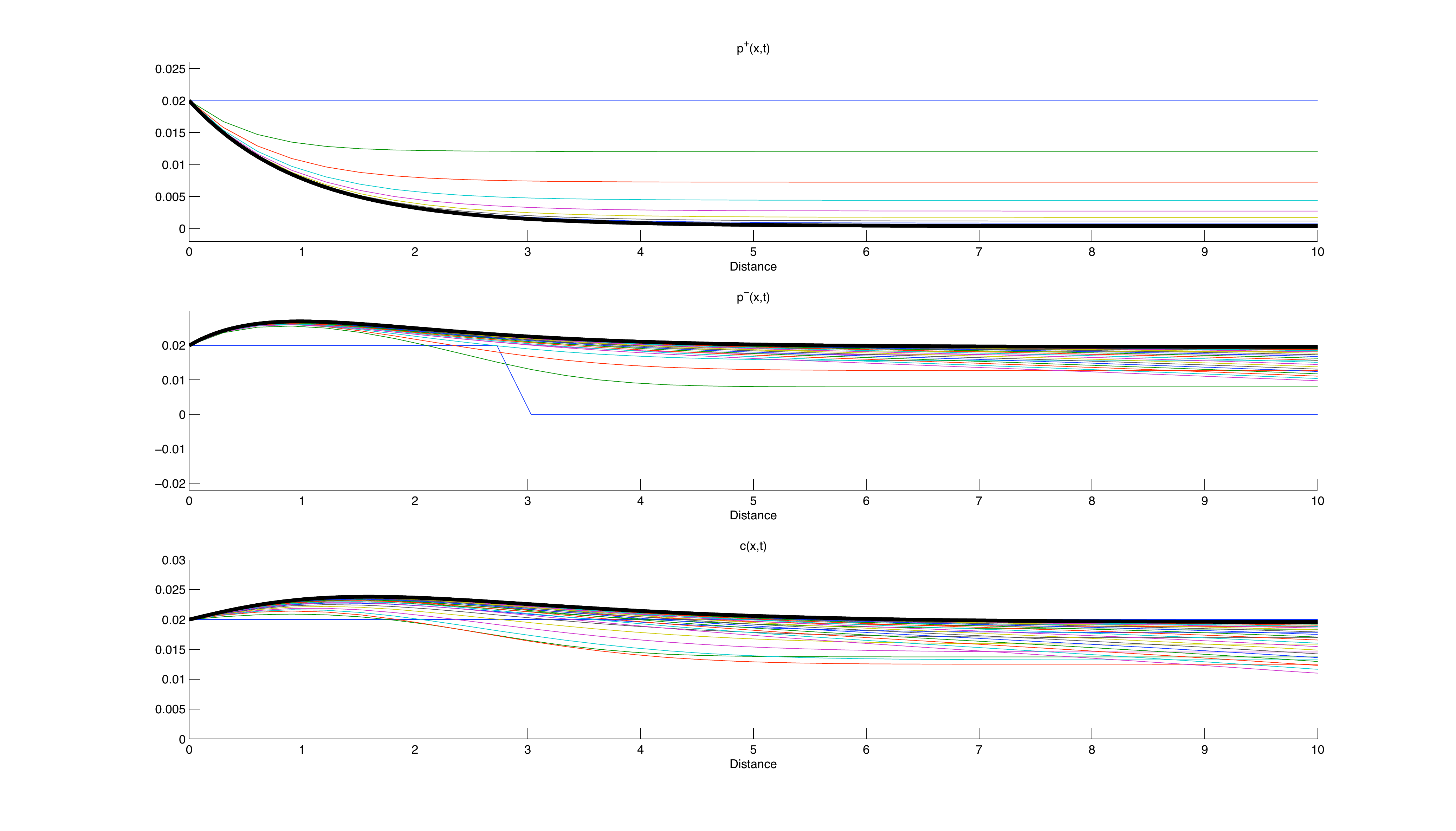}
\caption{Selection of endstate (iii) nonzero limiting concentration.
 }
\end{figure}

\section{Conclusion}
We have applied new mathematical tools to the question of how microtubule dynamic properties contribute to microtubule array construction.  Our results show the complex nature of the predicted end-states for a relatively simple model of microtubule dynamics including nucleation rate and tubulin concentration.  We emphasize that the mathematical tools introduced to handle he evident stochastic processes here are of general application, and should find use in related situations.  Originally developed to study stability of shock and boundary layers in gas dynamics, magnetohydrodynamics, and viscoelasticity
(See \cite{BHRZ,HLZ,CHNZ,HLyZ,BHZ,BLZ,BLeZ}),
these methods are capable of handling systems of essentially arbitrary complexity.
It is our hope that these methods and their application will open the way to the study of more complicated and realistic Biological models.
We view the present study {\it mainly as a feasibility study
for the application of new mathematical tools} transferring technology originating in the study of continuum mechanics to the area of biological modeling.

An important further direction from both Mathematical and Biological
point of view is the study of the singular perturbation limit as $d\to 0$.
In our numerics, we have mostly investigated the nonphysical case
$d\sim D$ conducive to good numerical conditioning,
whereas $d<<D$ in applications.
An important direction for further investigation, discussed in
preliminary fashion in Appendix \ref{singexist}, would be the rigorous
analysis of the singular limit $d\to 0$.
This should in principle be possible using the same body of techniques
applied here; see \cite{Z1} for an analysis in a similar spirit.
As discussed at the beginning of Section \ref{behavior}, another
very interesting open problem would be to establish stability of
steady state profiles with respect to perturbations decaying exponentially
only in the total density $p$ and not in other modes.
This should be accessible by the techniques of \cite{HZ,Z1}.


\appendix

\section{Resultant computations}\label{s:resultant}
\subsection{ Resultant Procedure for $c_+=0$}
\be
\ba
\\&
q_1(x)= \Big(Dd^2\Big)x^2+\Big(d^2u^+p^+_+ +d^2k+ f^+_-dD+ u^+c_+\nu^-D+d\omega c_+D\Big)x \\&
  +\Big(d\omega c_+k -d\omega p^-_+\nu^- +u^+c_+\nu^-k + d\omega c_+u^+p^+_+ +f_+^-dk +f^-_+du^+p^+_+-u^+c_+d\omega p^-_+\Big)\\&
q_2(x)= \Big(Dd\nu^--dDu^+c\Big)x^2-\Big(-f^+_-\nu^-D+ u^+c^2_+\omega D
+ u^+c_+dk- dv^-u^+p^+_+- d\nu^-k\Big)x \\ &+\Big(-f_+^-\nu^-k-u^+c_+^2\omega k\Big)
\ea
\ee

For $c_+=0$ these two reduce to
\be
\ba
\\&
q_1(x)= \Big(Dd^2\Big)x^2+\Big(d^2k+ f^+_-dD\Big)x
  +\Big(
f_+^-dk \Big)\\&
q_2(x)= \Big(Dd\nu^-\Big)x^2+\Big(f^+_-\nu^-D+ d\nu^-k\Big)x +\Big(-f_+^-\nu^-k\Big)
\ea
\ee
Subtracting $1/\nu^-$ times the second equation from $1/d$ times the
first yields a constant,
$$
2f^-_+k>0
$$
from which we may conclude that $q_1$ and $q_2$ have no common root.
{\it in particular, $q$ does not have any pure imaginary roots.}

\subsection{ Resultant Procedure for $c_+>0$}
Assuming $D=d=1$, $\nu^-=k=1$, $f^-_+=0$, $u^+=0$ and $\omega c_+= \omega p^+_+= \omega p^-_+=1$, we get
\be
\ba
\\&
q_1(x)= x^2+2x
 \\&
q_2(x)= x^2+x
\ea
\ee
which don't have a common positive root.

\section{The singular perturbation limit}\label{singexist}

In practice, it is important to take into account the singular perturbation
structure imposed on the problem by the small parameter $d$.
Otherwise, numerics will be destroyed when we try to go into the physical
parameter range $d<<1$, for one thing. For another, there is an advantage to singular limits, in that they tend to reduce to problem to composite problems that are sometimes explicitly solvable.
Finally, $d$ is often set to zero in applications, and it is important
to justify that this approximation is valid.
In this appendix, we initiate the discussion by a singular perturbation
study of the profile existence problem.  A parallel study of the
stability problem should be carried out, but appears to be more difficult.



First observe, if $d=0$, then the number of boundary conditions
should be only {\it two} at the boundary $x=0$: one for $c$,
which satisfies a parabolic equation, and one for the single
incoming mode among hyperbolic $p^\pm$.
So, we have too many boundary conditions for the slow problem $d=0$.

We may conclude from this that there is a ``true''
(i.e., in general nonconstant) boundary layer of thickness
$\sim d$  matching full parabolic to hyperbolic--parabolic boundary
conditions.  This satisfies the ``fast'' problem
obtained by rescaling $x\to \tilde x:= x/d$, and involves
second and first-order derivs. but {\it not} zero-order terms.
That is, the $p^\pm$ equations become in $\tilde x$ coordinates,
setting $d=0$, just
$$
p^+_{xx}= (\nu^+p^+)_x,
\quad
p^-_{xx}= -(\nu^-p^-)_x,
$$
a pair of conservation laws, and constant-coefficient as well.
Integrating up, we get
$$
p^+_x=\nu^+(p^+- p^+(+\infty)),
\quad
p^-_x=-\nu^-(P^-- -p^-(+\infty)),
$$
from which we find $p^+\equiv constant$, but $p_-(+\infty)$
is arbitrary, with an exponentially decaying solution.
The $c$ equation degenerates on the other hand to $c\equiv
{\rm constant}$.

So, the correct boundary data for the slow problem are
$p^+(0)$ prescribed (from original boundary conditions), $c$ prescribed (from
original boundary conditions, either Dirichlet or Neumann type preserved).

Then we solve the slow problem obtained by setting $d=0$ in
the original coordinates:
\begin{equation}\label{h1}
0=-\frac{\partial (\nu^+p^+(x,t))}{\partial x}-f^-_+p^+(x,t)+ f^+_-p^-(x,t)
\end{equation}
\begin{equation}
0=\nu^-\frac{\partial p^-(x,t)}{\partial x}+f^-_+p^+(x,t)- f^+_-p^-(x,t)
\end{equation}
\begin{equation}
0=-kc(x,t)+ \nu^- p^-(x,t)-\nu^+ p^+(x,t)+ D\frac{\partial^2c(x,t)}{\partial{x^2}}.
\end{equation}

Observing (by addition of the first two equations) that
\begin{equation}\label{relz}
\nu^+p^+-\nu^-p^-\equiv \alpha,
\end{equation}
$\alpha$ constant,
we may rewrite these equations after some rearrangement as
\begin{equation}\label{twoeq}
\begin{aligned}
(\nu^+p^+)_{x}
&=\Big( \frac{\omega c}{\nu^-  } -\frac{f^-_+ }{ u^+ c}\Big) (\nu^+p^+)
 - \frac{\omega c \alpha}{\nu^-} \\
\Big(c+ \frac{\alpha}{k}\Big)_{xx}&= \frac{k}{D}\Big(c + \frac{\alpha}{k}\Big).
\end{aligned}
\end{equation}
together with \eqref{relz}.
{\it Note that \eqref{relz} already selects the value $p^-(0)$
as a function of $p^+(0)$ and $c(0)$, in principle determining
completely the solution with no further computation,}
provided one exists.  The rest of the analysis consists in
verifying that \eqref{relz} is indeed consistent with existence.

The second, constant-coefficient equation together with the requirement
that $c$ be bounded as $x\to +\infty$ may be exactly solved as
$$
c(x) + \alpha/k= e^{-x\sqrt{k/D}}(c(0) + \alpha/k),
$$
$\alpha$ arbitrary, since only
the exponentially decaying mode is allowable.
In particular, we find that
\begin{equation}\label{calpha}
c(+\infty)=-\alpha/k.
\end{equation}
From \eqref{genass}, we find
that the coefficient
$\Big( \frac{\omega c^2u^+-f^-_+ \nu_-}{\nu^- u^+ c}\Big)$ of
$ (\nu^+p^+)$ in the first equation is negative as
$x\to +\infty$, so the stable manifold has a mode
in direction $p^+$.

Note that, by this observation,
the dimension of the phase space for the slow problem is four,
while the dimension of the stable manifold about a rest state is
two (one decaying direction for $c$,
and one decaying direction for $p^+$) and the dimension of the manifold
of solutions satisfying the (two) initial conditions is three
(two free parameters plus direction of spatial evolution).
Thus, as for the full problem, the dimension of the intersection is
generically a union of finitely many curves/steady state solutions;
in particular, we expect a connection to each possible endstate.

\medskip

We now investigate in detail the possible solutions.
There are two cases.

\medskip

{\it Case I. ($c(+\infty)=0$)}
In this case, also $\alpha=0$ by \eqref{calpha}, whence
$p^-(+\infty)=0$ by \eqref{relz} and $\nu^+(+\infty)=u^+c(+\infty)=0$.
Noting that in this case $c(x)= e^{-x\sqrt{k/D}}c(0)$,
we may rewrite the first equation of \eqref{twoeq}
as
$$
p^+_{x} = \frac{1}{u^+}
\Big( \frac{\omega c}{\nu^-} -
\frac{f^-_+ }{ u^+ c} + \sqrt{\frac{k}{D}}\Big) p^+
$$
to see that also $p^+\to 0$ as $x\to +\infty$,
as the solution of a homogeneous first-order
scalar equation with negative coefficient near $+\infty$.
Moreover, $(p^+,p^-)$ go to zero superexponentially, as
$e^{e^{-x\sqrt{k/D}}}$, while $c$ goes to zero only exponentially.
Such a slow solution exists for any initial data $p^+(0)$, $c^+(0)$,
hence a corresponding matched boundary-layer solution exists for
any Dirichlet data $p^+(0)$, $p^-(0)$, $c^+(0)$.

Neumann data: in this case, the only bounded solution of the $c$
equation satisfying the boundary conditions seems to be the
trivial solution $c\equiv 0$, from which \eqref{relz} gives
$p^-\equiv 0$.  But then \eqref{h1} yields $p^+\equiv 0$ as well,
i.e., only the trivial solution is possible, and we cannot generate
a slow solution of this type for general initial data.
\medskip

{\it Case II. ($c(+\infty)>0$)}
Again, only the decaying mode of $(c+\alpha/k)$ is possible, so that
the limit as $x\to +\infty$ of $c$ is $c(+\infty)=-\alpha/k$,
or
$$
\alpha=-kc(+\infty)<0.
$$
From \eqref{relz}, we obtain then $p^-(+\infty)=(c/\nu^-)(u^+p^+(+\infty)+k)$,
and the first equation of \eqref{twoeq} yields
at $+\infty$ the equilibrium condition
$p^+=\frac{\omega k c^2}{\nu^-f^-_+-u^+\omega c^2}$.
Solving for $p_-(+\infty)$ using \eqref{relz}, we obtain the full conditions
for $\cE$, as we should.
Likewise, by similar computations as above, we obtain a solution for
each initial data and endstate.

(Note: this is not a ``generic'' result as in the general case.
In other words, the assumption $d\to 0$ precludes nongeneric behavior.)

Neumann data: Again, the only bounded solution of the $c$
equation satisfying the boundary conditions is the
trivial solution $c\equiv -\alpha/k$,
greatly simplifying the analysis.
From \eqref{relz}, we obtain then $p^-=(c/\nu^-)(u^+p^++k)$,
and the first equation of \eqref{twoeq} becomes
$$
p^+_{x}
=\Big( \frac{\omega c^2u^+-f^-_+ \nu_-}{\nu^- u^+ c}\Big) p^+
 + \frac{\omega ck }{\nu^- u^+},
$$
yielding at $+\infty$ the equilibrium condition
$p^+=\frac{\omega k c^2}{\nu^-f^-_+-u^+\omega c^2}$.
Solving for $p_-(+\infty)$ using \eqref{relz}, we obtain the full conditions
for $\cE$, as we should.
Assuming \eqref{genass}, we again obtain a connection
for each endstate in $\cE$, as expected.
So, this case appears to permit profiles, also for Neumann data.

\begin{rems}
\textup{
1. Neumann boundary conditions do not seem consistent in case $\alpha=0$
corresponding to $c_+=0$.
Thus, this solution may (for Neumann conditions) be expected to
disappear in the singular limit $d\to 0$.
}

\textup{
2. At a heuristic level, one could just drop one boundary
condition and set $d=0$, it appears,
without much difference in profile behavior.
However, for stability studies one should not neglect the fast
(inner layer) dynamics near the boundary.
}
\end{rems}

In principle, stability should likewise be determinable
from the singular perturbation structure in the limit as $d\to 0$,
similarly as in \cite{Z1}. However, we do not attempt to carry
this out here.
Some such analysis should be carried out in justification of
the approximation $d=0$ commonly used in practice.

\end{document}